\documentclass[11pt]{article}
\usepackage[a4paper,margin=1in]{geometry}
\usepackage{amsmath,amssymb,amsthm,mathtools,booktabs,tabularx,array,enumitem,microtype,tikz}
\usetikzlibrary{arrows,positioning,shapes.geometric,calc}
\usepackage{xcolor,graphicx}
\usepackage{subfigure}
\numberwithin{figure}{section}
\numberwithin{table}{section} 
\usepackage[T1]{fontenc}\usepackage{lmodern}
\usepackage[colorlinks=true,linkcolor=blue,citecolor=blue,urlcolor=blue]{hyperref}
\numberwithin{equation}{section}
\newtheorem{theorem}{Theorem}[section]\newtheorem{lemma}[theorem]{Lemma}
\newtheorem{proposition}[theorem]{Proposition}\newtheorem{corollary}[theorem]{Corollary}
\theoremstyle{definition}
\newtheorem{assumption}[theorem]{Assumption}\newtheorem{definition}[theorem]{Definition}
\theoremstyle{remark}
\newtheorem{remark}[theorem]{Remark}
\newcommand{\R}{\mathbb R}\newcommand{\Pp}{\mathbb P}\newcommand{\E}{\mathbb E}
\newcommand{\supp}{\operatorname{supp}}
\newcommand{\dd}{\,\mathrm d}

\title{\textbf{Support of stochastic porous medium equations on bounded domains}}
\author{Guangying Lv$^{a}$ and Jinlong Wei$^b$\footnote{Corresponding author. weijinlong.hust@gmail.com (J. Wei)}\\  {\small \it $^a$ School of Mathematics and Statistics, Nanjing University of Information}\\ {\small \it  Science and Technology, Nanjing 210044, China}\\ {\small \it $^b$School of Statistics and Mathematics, Zhongnan University of}\\ {\small \it Economics and Law, Wuhan 430073, China}}
\date{}
\begin{document}\maketitle
\begin{abstract}
We study the spatial support of stochastic porous medium equations
on bounded domains. Under a local nondegeneracy condition, additive
noise gives full distributional support at every fixed positive time.
For the Dirichlet problem with nonlocal multiplicative noise, we prove
that the solution is either identically zero or has full support;
excluding the zero state requires a separate non-hitting assumption.
For a Neumann model, nonzero initial mass excludes the zero state and
yields full support. The proofs combine Galerkin approximation with
terminal Malliavin estimates. Numerical results illustrate the
propagation mechanisms.


\end{abstract}

\section{Introduction}
Let $\mathcal{O}\subset\R^d$ be bounded, connected, with $C^2$ boundary.  The
deterministic problem
\begin{equation}\label{1.1}\left\{\begin{array}{lll}
 \partial_t u(t,x)=\Delta(|u(t,x)|^{m-1}u(t,x)),\ \ & t>0,x\in\mathcal{O},\\
 u(t,x)=0,\ \  \ & x\in{\partial\mathcal{O}},\\
 u(0,x)=u_0,\ & x\in\mathcal{O},
 \end{array}\right.
\end{equation}
is degenerate at $u=0$ and has finite propagation when $m>1$. Its diffusion
coefficient vanishes at the vacuum, and this degeneracy produces a free
boundary and finite speed of propagation: bounded nonnegative initial data supported strictly inside the domain leave
an interior vacuum for sufficiently small positive times \cite{Abook}. Carrillo et al. \cite{CGT2004} used mass transportation techniques to give a simple proof of the finite speed of propagation of
the solution to the one-dimensional porous medium equation \eqref{1.1}.

\smallskip
A natural question is whether stochastic perturbations destroy this
localization phenomenon. More precisely, which noise structures can make the porous medium equation
(PME) lose finite propagation? The answer depends on how the noise acts
on the vacuum. We first recall the relevant results.

\smallskip
Barbu and R\"ockner \cite{BarbuRockner2012} studied stochastic porous medium equations with finite-dimensional linear multiplicative noise
\begin{equation}\label{1.2}\left\{\begin{array}{llll}
 \dd u(t,x)-\Delta(|u(t,x)|^{m-1}u(t,x))\dd t
 =
 \sum\limits_{k=1}^{N}\mu_k e_ku(t,x) \dd\beta_k(t), \ \  & t>0,x\in\mathcal{O},\\
 u(t,x)=0, \ \ & x\in \ \partial\mathcal{O},\\[0.1cm]
 u(0,x)=u_0(x), \ \ & x\in \ \mathcal{O},
   \end{array}\right.
\end{equation}
where $\beta_k$ is a sequence of independent Brownian motions on a filtered probability space $(\Omega,\mathcal{F},\mathcal{F}_t,\mathbb{P})$, $e_k$ is an orthonormal system in $L^2(\mathcal{O})$, and $\{\mu_k\}$ is a sequence of nonnegative numbers. Under their assumptions, including
$1<m\le5$ and $d\le3$, they proved almost-sure localization on a sufficiently
small random time interval.  Their proof used a random
exponential transformation that converts the stochastic partial differential equation (SPDE) into a degenerate
deterministic partial differential equation (PDE) with random coefficients, followed by a local energy method, also see \cite{BDPR2008,BDPR2009}. More precisely,
let
\[
 v(t,x)=e^{\mu(t)}u(t,x),
 \qquad
 \mu(t)=-\sum_{k=1}^N\mu_ke_k\beta_k(t),
\]
then the stochastic partial differential equation \eqref{1.2} becomes a deterministic PDE with random coefficients.

\smallskip
Gess \cite{Gess2013} proved finite speed of propagation for porous medium equations perturbed
by linear multiplicative space-time rough signals,
\begin{equation*}\left\{\begin{array}{llll}
 \dd u(t,x) = \Delta(|u(t,x)|^m\operatorname{sgn}u(t,x))\dd t
 + \sum_k f_ku(t,x)\circ\dd z^k(t),\ \ & t>0,x\in\mathcal{O},\\
 u(t,x)=0, \ \ & x\in\partial\mathcal{O},\\
 u(0,x)=u_0(x), \ \ & x\in \mathcal{O},
\end{array}\right.\end{equation*}
where rough driving signals $z^k\in C([0,T];\mathbb{R})$ and diffusion coefficients $f_k\in C^\infty(\mathcal{O})$. Moreover, Gess obtained the explicit and optimal estimates for the speed of
propagation. Hensel \cite{H2021} established finite time extinction with probability one for weak solutions of the Cauchy-Dirichlet problem for the one-dimensional stochastic porous medium equation with Stratonovich transport noise and compactly supported smooth initial datum. Fischer and Gr\"un developed a stochastic energy method for finite propagation and waiting-time
phenomena under linear multiplicative noise in dimensions up to three
\cite{FischerGrun2015}. They obtained waiting-time criteria and adapted the method to
finite speed of propagation. This gives further evidence that stochasticity
alone does not remove the free boundary when the noise remains local and
vacuum preserving. These results establish preservation of localization for the local coefficients and hypotheses considered in those works.

\smallskip
The common structural feature is
\begin{equation*}
 u=0
 \quad\Longrightarrow\quad
 f_ku=0,\ \ or\ \ \mu_ke_k u=0.
\end{equation*}
Thus the stochastic coefficient itself vanishes in a vacuum region.  The
random forcing may move the interface or alter its speed, but it does not
directly inject a disturbance into a region where the state vanishes.
Thus the noise is vacuum preserving.  These results are therefore not in
conflict with the proposed infinite-propagation mechanism.

\smallskip
On the other hand, Stan, del Teso and V\'azquez \cite{StanDelTesoVazquez2014,STV2016} studied the porous
medium equations with {\bf fractional pressure} or {\bf nonlocal pressure}
  \begin{equation*}\left\{\begin{array}{lll}
  \partial_tu(t,x)=\nabla\cdot(u^{m-1}(t,x)\nabla(-\Delta)^{-s}u(t,x)), \ \ & t>0,\\
  u(t,x)|_{t=0}=u_0(x), \ \ & x\in\mathbb{R}^d,
  \end{array}\right.\end{equation*}
where $m>1,0<s<1$. The main results of \cite{StanDelTesoVazquez2014,STV2016} show that spatial
nonlocality in the deterministic flux can change the finite speed
propagation regime. For $1<m<2$, the infinite-propagation result in
\cite[Theorem 1.4 (b)]{STV2016} is one-dimensional and is stated first as
a positive-measure tail property; pointwise positivity requires additional
symmetry and monotonicity assumptions. The case $m=1$ reduces to a fractional
heat equation. Further existence results are given in \cite{STV2019}.  Our equation keeps the classical local porous
medium drift and moves the nonlocal mechanism into the stochastic
state-to-noise operator.

\smallskip
We investigate two noise structures that can activate an interior vacuum.
One is \emph{additive noise}; the other is \emph{nonlocal multiplicative noise}. Additive noise generally does not preserve nonnegativity. We therefore
use distributional spatial support; when the state is a function, this agrees
with the support of its absolute value. Brownian nonlocal multiplicative
noise need not preserve nonnegativity either. For nonlocal multiplicative noise, we are motivated by \cite{StanDelTesoVazquez2014,STV2016}.
The main difference between local and nonlocal multiplicative noise can be interpreted as: for an open bounded domain $\mathcal{O}$,
\[
 B_k(u)(t,x):=g_k(x)u(t,x)
 \quad\Longrightarrow\quad
 u|_O=0\Rightarrow B_k(u)|_O=0
\]
for local multiplicative noise, where $O\Subset\mathcal O$ is a
nonempty open set, whereas
\[
 B(u)(t,x):=\int_{\mathcal{O}}K(x,y)u(t,y)\dd y
\]
may satisfy
\[
 u|_O=0,\quad u\not\equiv0
 \quad\Longrightarrow\quad
 B(u)|_O\ne0.
\]
For a strictly positive convolution kernel and a nonnegative nontrivial
state,
\[
 \supp u\Subset\mathcal{O}
 \quad\Longrightarrow\quad
 \supp_{\mathcal O} B(u)=\overline{\mathcal O}.
\]
This is the mechanism absent from the known local multiplicative
finite speed propagation theory.

\smallskip
The above discussions motivate the following questions:

\smallskip
\noindent\textbf{Question 1:}
\emph{Can additive noise destroy the compact support in the following sense, for every fixed $t>0$,
\[
 \Pp(\supp u(t,\cdot)=\overline{\mathcal O})=1?
\]}

\noindent\textbf{Question 2:}
\emph{Can one construct a multiplicative stochastic perturbation for which
the deterministic porous medium equation has finite speed of propagation,
whereas the stochastic equation has instantaneous full spatial support?}

\smallskip
The central idea of this paper is that the relevant distinction is not merely
between deterministic and stochastic equations, nor between finite- and
infinite-dimensional Wiener processes. Rather, one should distinguish
between \emph{vacuum-preserving} and \emph{vacuum-destroying} noise
coefficients. Theorem~\ref{thm2.16} makes this mechanism rigorous for additive noise.
For multiplicative noise, we propose the genuinely nonlocal multiplicative model
\begin{equation}\label{1.3}
 \dd u(t,x)
 =
 \Delta\Phi(u(t,x))\dd t+(K*u)(t,x)\dd\beta_t,
 \qquad
 \Phi(r)=|r|^{m-1}r,
 \qquad m>1,
\end{equation}
where
\[
 (K*u)(t,x)=\int_{\mathcal{O}}K(x-y)u(t,y)\dd y.
\]
If $K>0$ almost everywhere, and $u\ge0$ and $u\not\equiv0$, then
\[
 (K*u)(t,x)>0
 \quad\text{for a.e. }x\in\mathcal{O},
\]
even when $u(t,x)=0$ on a large open set in spatial variable.  Hence the stochastic coefficient
directly sees every vacuum region. This is stronger than merely choosing a nonlocal covariance operator.  For example, $u Q^{1/2}\dd W$ may have a spatially nonlocal covariance, but the pointwise factor $u(x)$
still forces the noise to vanish at every vacuum point.  In
\eqref{1.3}, by contrast, the \emph{state-to-noise map} itself is
nonlocal.

\smallskip
The propagation mechanisms and the conclusions proved here are summarized in Table~\ref{tab1.1}.
\begin{table}[htbp]
\centering
\small
\renewcommand{\arraystretch}{1.25}
\begin{tabularx}{\textwidth}{
 >{\raggedright\arraybackslash}p{0.29\textwidth}
 >{\raggedright\arraybackslash}p{0.10\textwidth}
 >{\raggedright\arraybackslash}p{0.19\textwidth}
 >{\raggedright\arraybackslash}X}
\toprule
\textbf{Model}
& \textbf{Noise on vacuum}
& \textbf{Spatial mechanism}
& \textbf{Propagation behavior} \\
\midrule
$\partial_tu=\Delta u^m$
& no noise
& local degenerate diffusion
& finite speed for $m>1$ \cite{Abook} \\
$\partial_tu=\nabla\cdot(u^{m-1}\nabla(-\Delta)^{-s}u)$
& no noise
& nonlocal degenerate diffusion
& infinite propagation for $1<m<2$ in dimension one under the stated hypotheses \cite{STV2016} \\

$\dd u=\Delta\Phi(u)\dd t+u \dd W$
& zero
& local multiplicative
& finite propagation under suitable hypotheses \cite{BarbuRockner2012} \\

$\dd u=\Delta\Phi(u)\dd t+\sum_k f_ku\circ\dd z_k$
& zero
& local multiplicative rough noise
& finite propagation \cite{Gess2013} \\

$\dd u=\Delta\Phi(u)\dd t+Q^{1/2}\dd W$
& generally nonzero
& \textbf{additive, spatially correlated}
& \textbf{full support under local nondegeneracy; positivity is generally lost} \\

$\dd u=\Delta\Phi(u)\dd t+(K*u)\dd\beta_t$
& generally nonzero
& \textbf{nonlocal multiplicative}
& \textbf{vacuum activation; full support needs further hypotheses} \\
\bottomrule
\end{tabularx}
\caption{Local and nonlocal propagation mechanisms. }
\label{tab1.1}
\end{table}

The main contributions, in this paper, are

\smallskip
1) The full-domain support is obtained for additive noise. The additive-noise Malliavin method is used and comparing with \cite{BakhtinMattingly}, our differentiation step instead applies a finite-dimensional monotone-SDE theorem before taking a weak Galerkin limit.

\smallskip
2) Dirichlet full-support or zero-state dichotomy is given for Dirichlet boundary problem. In order to get the Dirichlet full-support, the backward non-hitting property should be added.

\smallskip
3) Different from Dirichlet boundary problem, the full support for the Neumann mass model with nonzero initial mass is obtained.

\smallskip
The rest of the paper is organized as follows. Section~2 proves the additive result and develops the common approximation
framework. Section~3 establishes the mixed-norm terminal estimate and the
Dirichlet dichotomy. Section~4 proves the Neumann result by conditioning on
the mass noise. A simple numerical example is given is Section 5. Conclusion is given in Section~6.

\section{Additive noise and the approximation framework}

Throughout the paper, $m>1$, $p=m+1$, and $\mathcal O\subset\mathbb R^d$
is bounded and connected with $C^2$ boundary. All initial conditions in the
Malliavin statements are deterministic. We write
$\langle v,\varphi\rangle$ for distributional duality and $\|v\|_q$ for
the $L^q(\mathcal O)$ norm. We first introduce the definition of spatial support.
\begin{definition}\label{def2.1}
For $v\in\mathcal D'(\mathcal O)$, let
\[
 S_{\mathcal O}(v)=\mathcal O\setminus\bigcup\{O\subset\mathcal O:
 O\text{ is open and }v|_O=0\}.
\]
Define
\[
\supp_{\mathcal O}v=\overline{S_{\mathcal O}(v)}^{ \mathbb R^d}.
\]
In particular, $\supp_{\mathcal O}0=\varnothing$, where $0$ means the zero function. Full-domain support means
$\supp_{\mathcal O}v=\overline{\mathcal O}$.
\end{definition}

Every support assertion below uses Definition~\ref{def2.1}. If $v$
is locally integrable, it vanishes as a distribution on an open set exactly
when it vanishes almost everywhere there. Thus this convention agrees with
the closed support of $|v|$, and allows signed states without imposing a
positivity property on the stochastic equation.

\smallskip
In this paper, we want to show the full-support of stochastic porous medium equation. In order that, we need the instantaneous full-domain propagation for any prescribed time.

\begin{definition}\label{def2.2}
A stochastic solution has instantaneous full-domain propagation at fixed
positive times if, for every fixed $t>0$,
\[
\Pp(\supp_{\mathcal O}u(t,\cdot)=\overline{\mathcal O})=1.
\]
\end{definition}

Definition~\ref{def2.2} concerns each prescribed time, so it does not
assert a single probability-one event valid for all positive real times.
It is also compatible with the homogeneous Dirichlet boundary condition:
a function can have zero boundary trace while being nonzero on every
interior open set. No pointwise strict positivity is included in the definition. We recall the deterministic benchmark before introducing noise. Let
$B_R(x_0)$ be the open Euclidean ball with radius $R$ centred $x_0$.

\begin{theorem}[Persistence of a deterministic interior vacuum]
\label{thm2.3}
Let $0\le v_0\in L^\infty(\mathcal O)$, and let $v$ be the nonnegative
weak solution of \eqref{1.1}. Suppose $v_0=0$ on
$B_R(x_0)\Subset\mathcal O$. For every $0<R_1<R$, there is $T_*>0$
such that $v(t)=0$ on $B_{R_1}(x_0)$ for $0<t<T_*$.
In particular, bounded nonnegative initial data supported strictly inside
$\mathcal O$ do not produce instantaneous full-domain propagation.
\end{theorem}

\begin{proof} The result is known and we only give the outline of proof.
The classical finite-propagation, or hole-filling, estimate for the porous
medium equation bounds the advance of a free boundary into an initial
vacuum by a function $\eta(t)$ tending to zero as $t\downarrow0$.
It follows either from the local energy method or from comparison with
porous-medium barriers; see \cite{Abook,Vazquez2007}. Applied to the
initially empty ball, it gives $v(t)=0$ on
$B_{R-\eta(t)}(x_0)$ as long as this radius is positive.
The bound depends on the initial $L^\infty$ norm and the local geometry,
and the zero Dirichlet condition supplies no incoming boundary source.

\smallskip
Choose $T_*$ so that $\eta(t)<R-R_1$ for $0<t<T_*$. Then the smaller
ball remains empty throughout that interval. If
$\supp_{\mathcal O}v_0\Subset\mathcal O$, its complement contains such
an interior ball, so the first assertion applies. This argument concerns
small positive times; it does not say that the support remains strictly
inside the domain for every finite time.
\end{proof}

The contrast with Theorem~\ref{thm2.3} will be that sufficiently
nondegenerate additive noise removes every interior vacuum at any fixed
positive time, including when the initial datum is compactly supported.

\subsection{Additive noise and well-posedness}
Let
\[
H_D=H^{-1}(\mathcal O)=(H_0^1(\mathcal O))^* \quad {\rm with}
\quad \|v\|_{H_D}^2=\langle v,(-\Delta_D)^{-1}v\rangle.
\]
The variational
triple is $V=L^p(\mathcal O)\subset H_D\subset V^*$. Define the positive
drift $A:V\to V^*$ by
\begin{equation}\label{2.1}
 {}_{V^*}\langle A(v),w\rangle_V=\int_{\mathcal O}\Phi(v(x))w(x) \dd x,
 \qquad \Phi(r)=|r|^{m-1}r.
\end{equation}
Thus $-A$ is the variational realization of $\Delta\Phi$ with homogeneous
Dirichlet boundary conditions. An $H_D$-continuous solution need not be an
$L^p$ function at every fixed time; this is one reason to use distributional
support.

\smallskip
We establish the estimates for affine noise, so that they can be reused.
Let $W$ be cylindrical on a separable Hilbert space $\mathcal U$, let
$g\in C([0,T];\mathcal L_2(\mathcal U,H_D))$ be deterministic, and let
$C\in\mathcal L(H_D,\mathcal L_2(\mathcal U,H_D))$. Write
$C_kv=(Cv)e_k$ for an orthonormal basis $(e_k)$ of $\mathcal U$, and assume
\begin{equation}\label{2.2}
 \sum_{k\ge1}\|C_kv\|_{H_D}^2\le L\|v\|_{H_D}^2,
 \qquad v\in H_D.
\end{equation}
Set $\Sigma(s,v)=g(s)+Cv$ and consider
\begin{equation}\label{2.3}
 \dd X(s)+A(X)\dd s=\Sigma(s,X)\dd W_s,\qquad X(0)=x\in H_D.
\end{equation}
Here we let $u$ replaced by $X$ in order to consider the abstract SDEs. The diffusion term $\Sigma(s,X)$ contains the additive noise and linear multiplicative noise.
For a smooth test $\varphi$, the notation $\Sigma(s,v)^*\varphi\in\mathcal U$
means the vector with coordinates $\langle\Sigma(s,v)e_k,\varphi\rangle$.
It is defined by duality with $H_0^1$, rather than by identifying that test
with an element of $H_D$. If the noise takes values in $L^2$, this agrees
with the usual $L^2$ adjoint notation.

\smallskip
The additive equation is the special case $g(s)=G$ and $C=0$, where
$G\in\mathcal L_2(\mathcal U,H_D)$ is deterministic.

\begin{assumption}[Additive local nondegeneracy]\label{ass2.4}
For every nonzero $\varphi\in C_c^\infty(\mathcal O)$,
$\|G^*\varphi\|_{\mathcal U}>0$.
\end{assumption}

Assumption~\ref{ass2.4} asks that every localized scalar observation
see at least one noise direction. It does not require a uniform lower
bound over all tests. The following familiar trace-class example shows
that the assumption is compatible with spatially correlated noise.

\begin{proposition}[A spectral sufficient condition]\label{prop2.5}
Let $(e_k)$ be a complete orthonormal basis of $L^2(\mathcal O)$ and let
$Qe_k=q_ke_k$, where $q_k>0$ and $\sum_kq_k<\infty$. Then
\[
G=Q^{1/2}:L^2(\mathcal O)\to L^2(\mathcal O),
\]
also viewed as an
$H_D$-valued operator, satisfies Assumption~\ref{ass2.4}.
\end{proposition}

\begin{proof}
The Hilbert--Schmidt norm into $L^2$ is
\[
\|G\|_{\mathcal L_2}^2=\sum_kq_k<\infty.
\]
The continuous embedding
$L^2\subset H_D$ also makes $G$ Hilbert--Schmidt into $H_D$.
For a nonzero test $\varphi$, completeness supplies a $k_0$ with
$(\varphi,e_{k_0})_{L^2}\ne0$. Consequently
\[
 \|G^*\varphi\|_{\mathcal U}^2
 =\sum_k q_k|(\varphi,e_k)_{L^2}|^2
 \ge q_{k_0}|(\varphi,e_{k_0})_{L^2}|^2>0.
\]
Here the adjoint is used in the duality sense specified above, which gives
exactly the displayed coordinates.
\end{proof}

For $F_{t,\varphi}=\langle u(t),\varphi\rangle$, the scalar weak equation is
\begin{equation}\label{2.4}
 F_{t,\varphi}=\langle u_0,\varphi\rangle
 +\int_0^t\!\int_{\mathcal O}\Phi(u(s))\Delta\varphi \dd x\dd s
 +\int_0^t(G^*\varphi,\dd W_s)_{\mathcal U}.
\end{equation}
Applying It\^o's isometry to the deterministic integrand $G^*\varphi$, we get
\[
 \E\left|\int_0^t(G^*\varphi,\dd W_s)_{\mathcal U}\right|^2
 =t\|G^*\varphi\|_{\mathcal U}^2>0.
\]
The above equality explains why the noise can act inside a vacuum:
its scalar coefficient is independent of the state. It does not yet prove
$\Pp(F_{t,\varphi}=0)=0$. The drift in \eqref{2.4} depends
on the same Wiener path, so the observable is not an independent random
shift of the displayed Gaussian term. The approximation and Malliavin
arguments below establish the needed nondegeneracy of the complete
observable.

\smallskip
To establish well-posedness, we use a Galerkin approximation based on
the eigenfunctions of the Dirichlet Laplacian. Let
$(\rho_j,\lambda_j)_{j\geq1}$ satisfy
\[
-\Delta\rho_j=\lambda_j\rho_j \quad\text{in }\mathcal O,
\qquad
\rho_j=0 \quad\text{on }\partial\mathcal O,
\]
where $(\rho_j)_{j\geq1}$ is an orthonormal basis of
$L^2(\mathcal O)$ and the eigenvalues are listed in nondecreasing order. Let $H_n=\operatorname{span}\{\rho_1,\ldots,\rho_n\}$ and let $P_n$ be the spectral projection. It is also orthogonal in $H_D$ and contractive on $H_0^1$. Let $\Pi_n$ project $\mathcal U$ onto its first $n$ coordinates. In this paper, we do not approximate the noise \cite{MZ2019} and only discretized the space.

\begin{lemma}[Projected test functions]\label{lem2.7}
For $\varphi\in C_c^\infty(\mathcal O)$, set $\varphi_n=P_n\varphi$.
Then $\varphi_n\to\varphi$ in $H_0^1$ and
\begin{equation}\label{2.5}
 \sup_n\|\Delta\varphi_n\|_\infty<\infty,
 \qquad \Delta\varphi_n\longrightarrow\Delta\varphi
 \quad\hbox{uniformly on }\mathcal O.
\end{equation}
The corresponding statement holds for the positive Neumann modes, with
$P_n\varphi\to\varphi-\overline\varphi$ in $H^1$.
\end{lemma}

\begin{proof}
A compactly supported smooth test belongs to the domain of every integer
power of the Dirichlet Laplacian. For the Neumann Laplacian this is true of
$\varphi-\overline\varphi$. The standard heat-kernel and eigenvalue estimates
on a bounded $C^2$ domain give
\[
\|\rho_j\|_\infty\le C(1+\lambda_j)^{\frac{d}{4}}
\]
and
$\lambda_j\ge c j^{2/d}$ for the positive modes. If $a_j=(\varphi,\rho_j)_{L^2}$,
then
\[
\sum_j(1+\lambda_j)^{2q}|a_j|^2<\infty
\]
for every integer $q$.
For $q>1+d/2$, Cauchy--Schwarz gives
\[
\sum_j\lambda_j|a_j|\|\rho_j\|_\infty<\infty.
\]
More explicitly, the
second factor in that inequality is bounded by
\[
 C\Big(\sum_j(1+\lambda_j)^{2+d/2-2q}\Big)^{\frac12}<\infty,
\]
while the first is
$(\sum_j(1+\lambda_j)^{2q}|a_j|^2)^{1/2}$.
The series for $\Delta\varphi_n$ therefore converges absolutely and
uniformly, with limit $\Delta\varphi$. Spectral convergence in the energy
norm proves the other assertions.
\end{proof}

Although $\varphi$ is compactly supported, $\varphi_n$ generally is not.
We will therefore estimate the linearized drift over the whole domain
$\mathcal O$. The uniform Laplacian bound in
Lemma~\ref{lem2.7} is what makes that estimate independent
of the Galerkin dimension.

\smallskip
Define $A_n:H_n\to H_n$ by
\[
(A_n(v),w)_{H_D}=\int_{\mathcal O}\Phi(v)w \dd x \quad {\rm for} \quad w\in H_n.
\]
We use the actual porous-medium nonlinearity in the finite-dimensional SDE
\begin{equation}\label{2.6}
 \dd X_n(s)+A_n(X_n)\dd s
 =P_n\Sigma(s,X_n)\Pi_n\dd W_s,\qquad X_n(0)=P_nx.
\end{equation}
No globally bounded derivative is imposed on $\Phi$. Since $m>1$, $\Phi$
is $C^1$, which suffices for the finite-dimensional differentiability result
used below. In particular, no incompatible combination of bounded derivative
and superquadratic coercivity enters the approximation.

\begin{proposition}[Well-posedness]
\label{prop2.8}
Under \eqref{2.2}, equation \eqref{2.3} has a
unique solution $X$ in $L^2(\Omega;C([0,T];H_D))\cap
L^p(\Omega\times(0,T)\times\mathcal O)$. The solutions of
\eqref{2.6} satisfy
\begin{equation}\label{2.7}
 \sup_n\left\{\E\sup_{s\le T}\|X_n(s)\|_{H_D}^2
       +\E\int_0^T\|X_n(s)\|_p^p \dd s\right\}<\infty.
\end{equation}
They converge weakly to $X$ in $L^p(\Omega\times(0,T)\times\mathcal O)\cap L^2(\Omega\times(0,T);H_D)$. For every fixed $t\in[0,T]$ and
$\varphi\in C_c^\infty(\mathcal O)$,
\begin{equation}\label{2.8}
 \langle X_n(t),\varphi\rangle\rightharpoonup
 \langle X(t),\varphi\rangle\quad\hbox{in }L^2(\Omega).
\end{equation}
\end{proposition}

\begin{proof} The proof is divided into three steps.

\smallskip
\emph{Step 1: the variational hypotheses and uniform energy.}
For \eqref{2.1}, semi-continuity follows from continuity of
$\Phi$, while
\[
 \langle A(v)-A(w),v-w\rangle
 =\int_{\mathcal O}[\Phi(v)-\Phi(w)](v-w) \dd x\ge0,
 \qquad \|A(v)\|_{V^*}\le\|v\|_p^{p-1}.
\]
Moreover $\langle A(v),v\rangle=\|v\|_p^p$ and
\[
\|\Sigma(s,v)\|_{\mathcal L_2}^2\le2\|g(s)\|_{\mathcal L_2}^2+2L\|v\|_{H_D}^2.
\]
These are the semi-continuity, monotonicity, coercivity and growth hypotheses
of the variational theorem \cite{LiuRockner2010}. In particular, we use its
polynomial variational setting, without attributing arbitrary $H_D$-valued
noise to a theorem that imposes spatial smoothing on the noise.

\smallskip
The Galerkin drift is $C^1$ and dissipative, and the noise is affine and
globally Lipschitz. It\^o's formula gives exactly
\[
 \begin{aligned}
 \|X_n(t)\|_{H_D}^2+2\int_0^t\|X_n(s)\|_p^p \dd s
 &=\|P_nx\|_{H_D}^2+\int_0^t\|P_n\Sigma(s,X_n)\Pi_n\|_{\mathcal L_2}^2 \dd s\\
 &\quad+2\int_0^t(X_n,P_n\Sigma(s,X_n)\Pi_n\dd W_s)_{H_D}.
 \end{aligned}
\]
Note that the uniformity in $n$ is used repeatedly below. Set
\[
 \Sigma_n(s)=P_n\Sigma(s,X_n(s))\Pi_n,\qquad
 M_n(t)=2\int_0^t(X_n(s),\Sigma_n(s)\dd W_s)_{H_D}.
\]
For $R>0$, let
\[
 \tau_R=\inf\{s\ge0:\|X_n(s)\|_{H_D}\ge R\}\wedge T .
\]
Taking the supremum in the preceding It\^o identity up to
$t\wedge\tau_R$ and then expectation gives
\[
 \begin{aligned}
 &\E\sup_{\theta\le t\wedge\tau_R}\|X_n(\theta)\|_{H_D}^2
 +2\E\int_0^{t\wedge\tau_R}\|X_n(s)\|_p^p\dd s\\
 &\le
 \|P_nx\|_{H_D}^2
 +\E\int_0^{t\wedge\tau_R}\|\Sigma_n(s)\|_{\mathcal L_2}^2\dd s
 +\E\sup_{\theta\le t\wedge\tau_R}|M_n(\theta)|.
 \end{aligned}
\]
Since $P_n$ and $\Pi_n$ are contractions,
\[
 \|\Sigma_n(s)\|_{\mathcal L_2}^2
 \le \|\Sigma(s,X_n(s))\|_{\mathcal L_2}^2
 \le 2\|g(s)\|_{\mathcal L_2}^2+2L\|X_n(s)\|_{H_D}^2.
\]
Moreover, the BDG inequality followed by Young's inequality yields
\[
 \begin{aligned}
 \E\sup_{\theta\le t\wedge\tau_R}|M_n(\theta)|
 &\le C\E\!\left[
   \sup_{\theta\le t\wedge\tau_R}\|X_n(\theta)\|_{H_D}
   \left(\int_0^{t\wedge\tau_R}
       \|\Sigma_n(s)\|_{\mathcal L_2}^2\dd s\right)^{1/2}\right]\\
 &\le \frac12\E\sup_{\theta\le t\wedge\tau_R}
       \|X_n(\theta)\|_{H_D}^2
   +C\E\int_0^{t\wedge\tau_R}
       \|\Sigma_n(s)\|_{\mathcal L_2}^2\dd s .
 \end{aligned}
\]
Consequently,
\[
 \begin{aligned}
 &\frac12\E\sup_{\theta\le t\wedge\tau_R}\|X_n(\theta)\|_{H_D}^2
 +2\E\int_0^{t\wedge\tau_R}\|X_n(s)\|_p^p\dd s\\
 &\le
 \|x\|_{H_D}^2
 +C\int_0^t\|g(s)\|_{\mathcal L_2}^2\dd s
 +CL\int_0^t
     \E\sup_{\theta\le s\wedge\tau_R}\|X_n(\theta)\|_{H_D}^2\dd s .
 \end{aligned}
\]
Gronwall's inequality gives
\[
 \sup_n\E\sup_{\theta\le T\wedge\tau_R}\|X_n(\theta)\|_{H_D}^2
 \le C_T\left(\|x\|_{H_D}^2+
       \int_0^T\|g(s)\|_{\mathcal L_2}^2\dd s\right),
\]
with a constant independent of $n$ and $R$.  Substitution into the
preceding energy inequality gives the same uniform control of
$\E\int_0^{T\wedge\tau_R}\|X_n(s)\|_p^p\dd s$.  Letting
$R\uparrow\infty$ and using Fatou's lemma for the supremum term and
monotone convergence for the nonnegative integral proves \eqref{2.7};
in particular, the stopped solutions cannot explode.

\smallskip
For later use we also record why higher $H_D$ moments are available.
Applying It\^o's formula to
\[
\Psi_q(y)=(1+\|y\|_{H_D}^2)^{q/2},\quad q\ge2,
\]
the dissipative drift
contributes a nonpositive term, while the first- and second-order
diffusion terms are bounded, by the linear-growth estimate above, by
$C_q(1+\|g(s)\|_{\mathcal L_2}^q+\Psi_q(X_n(s)))$.  The same stopped
BDG--Young--Gronwall argument therefore gives
\[
 \sup_n\E\sup_{s\le T}(1+\|X_n(s)\|_{H_D}^2)^{q/2}<\infty
\]
for every fixed $q\ge2$.

\smallskip
\emph{Step 2: identification of the weak Galerkin limit.}
The energy bound gives a weakly convergent subsequence of $X_n$ in $L^p$
and in $L^2(\Omega\times(0,T);H_D)$. The drifts are bounded in
$L^{p/(p-1)}(\Omega\times(0,T);V^*)$, so denote their weak limit by
$\overline A$. Affine linearity identifies the weak noise limit with
$\Sigma(\cdot,X)$: in a pairing with a fixed Hilbert--Schmidt test,
move $P_n$ and $\Pi_n$ to that test and use their strong convergence.
Testing the Galerkin equations against each fixed finite span gives
\[
 X(t)=x-\int_0^t\overline A(s) \dd s
             +\int_0^t\Sigma(s,X(s)) \dd W_s
 \quad\hbox{in }V^*.
\]
By using the variational It\^o formula, we can get an $H_D$-continuous version.
We next justify the fixed-time weak convergence needed in the Minty
argument, without assuming any path-space strong convergence.
The estimate \eqref{2.7} makes $\{X_n(t)\}_n$ bounded in
$L^2(\Omega;H_D)$ for every fixed $t$.  Hence a subsequence converges
weakly there to some $Y(t)$.  Let $h$ belong to a fixed finite
eigenspace.  For all sufficiently large $n$, testing the Galerkin
equation with $h$ gives
\[
 (X_n(t),h)_{H_D}
 =(P_nx,h)_{H_D}
 -\int_0^t {}_{V^*}\!\langle A(X_n(s)),h\rangle_V\dd s
 +\left(\int_0^t\Sigma_n(s)\dd W_s,h\right)_{H_D}.
\]
The time-integration map sends
$L^{p'}(\Omega\times(0,t);V^*)$ continuously into the corresponding
scalar $L^{p'}(\Omega)$ pairing with $h$, while stochastic integration is
a bounded linear map from the predictable
$L^2(\Omega\times(0,t);\mathcal L_2(\mathcal U,H_D))$ space into
$L^2(\Omega;H_D)$.  Thus, after multiplying by an arbitrary bounded
random variable $Z$, the weak convergences of $A(X_n)$ and $\Sigma_n$
pass through the two integrals.  The limit of the scalar pairing is
\[
 (x,h)_{H_D}
 -\int_0^t {}_{V^*}\!\langle\overline A(s),h\rangle_V\dd s
 +\left(\int_0^t\Sigma(s,X(s))\dd W_s,h\right)_{H_D}
 =(X(t),h)_{H_D}.
\]
Finite eigenspaces are dense in $H_D$, bounded random variables are
dense in $L^2(\Omega)$, and the variables $X_n(t)$ are uniformly bounded
in $L^2(\Omega;H_D)$.  Finite sums of simple tensors $Zh$ are therefore
dense in $L^2(\Omega;H_D)$, so $Y(t)=X(t)$.  Consequently
\[
 X_n(t)\rightharpoonup X(t)
 \qquad\hbox{in }L^2(\Omega;H_D)
\]
for every fixed $t\in[0,T]$.

\smallskip
We spell out the monotonicity argument because weak convergence by itself
does not identify a nonlinear drift. Choose $\kappa>L$. The weighted
Galerkin energy identity and contractivity of the noise projections give
\[
 \begin{split}
 J_n:={}&\E\int_0^T e^{-\kappa s}
 \left(2\langle A(X_n),X_n\rangle+\kappa\|X_n\|_{H_D}^2
                         -\|\Sigma(s,X_n)\|_{\mathcal L_2}^2\right)\dd s\\
 \le{}&\|P_nx\|_{H_D}^2-\E e^{-\kappa T}\|X_n(T)\|_{H_D}^2.
 \end{split}
\]
Lower semicontinuity of the terminal norm and the energy identity for the
candidate show that
\[
 \limsup_n J_n\le
 \E\int_0^T e^{-\kappa s}
 \left(2\langle\overline A,X\rangle+\kappa\|X\|_{H_D}^2
                              -\|\Sigma(s,X)\|_{\mathcal L_2}^2\right)\dd s.
\]
For a predictable $v\in L^p(\Omega\times(0,T);V)$, monotonicity and
the Lipschitz noise bound imply
\[
 \begin{split}
 0\le \E\int_0^T e^{-\kappa s}\bigl(
2\langle A(X_n)-A(v),X_n-v\rangle+\kappa\|X_n-v\|_{H_D}^2-\|\Sigma(s,X_n)-\Sigma(s,v)\|_{\mathcal L_2}^2\bigr)\dd s.
 \end{split}
\]
Expand this expression. Its quadratic part is $J_n$, and all remaining
terms pass to the limit by weak convergence. The preceding limsup bound
therefore yields
\[
 \begin{split}
 0\le \E\int_0^T e^{-\kappa s}\bigl(
2\langle\overline A-A(v),X-v\rangle+\kappa\|X-v\|_{H_D}^2-\|\Sigma(s,X)-\Sigma(s,v)\|_{\mathcal L_2}^2\bigr)\dd s.
 \end{split}
\]
Take $v=X-\varepsilon w$ and $v=X+\varepsilon w$, first with bounded
predictable $V$-valued $w$. Divide by $\varepsilon$ and let
$\varepsilon\downarrow0$. Hemicontinuity and the polynomial growth bound
give $\overline A=A(X)$ by the two opposite inequalities and density.
This is the required Minty identification.

\smallskip
For two solutions driven by the same noise and with the same initial
condition, It\^o's formula and monotonicity give
\[
\E\|X(t)-\widetilde X(t)\|_{H_D}^2
\le L\int_0^t\E\|X(s)-\widetilde X(s)\|_{H_D}^2 \dd s.
\]
Gronwall proves uniqueness. Thus every weakly convergent subsequence has
the same limit, and the whole Galerkin sequence converges in the asserted
weak topologies.

\smallskip
\emph{Step 3: the scalar observation at a prescribed time.}
One still needs a fixed-time assertion for Malliavin calculus. Let
$F_n(s)=\langle X_n(s),\varphi\rangle$. Its scalar equation has drift
\[
\int_{\mathcal O}\Phi(X_n)\Delta\varphi_n \dd x
\]
and noise
$\Pi_n\Sigma(s,X_n)^*\varphi_n$. By Lemma~\ref{lem2.7},
H\"older's inequality and \eqref{2.7},
\begin{equation}\label{2.9}
 \sup_n\E|F_n(t)-F_n(s)|
 \le C_\varphi\bigl(|t-s|^{\frac1p}+|t-s|^{\frac12}\bigr).
\end{equation}
For completeness, we indicate the estimate behind \eqref{2.9}.  The
drift increment satisfies, uniformly in $n$,
\[
 \begin{aligned}
 \E\left|\int_s^t\!\!\int_{\mathcal O}
      \Phi(X_n(q))\Delta\varphi_n\dd x\dd q\right|
 &\le C_\varphi
 \E\int_s^t\|X_n(q)\|_p^{p-1}\dd q\\
 &\le C_\varphi |t-s|^{1/p}
 \left(\E\int_s^t\|X_n(q)\|_p^p\dd q\right)^{(p-1)/p}
 \le C_\varphi |t-s|^{1/p},
 \end{aligned}
\]
where \eqref{2.5} and \eqref{2.7} were used.  The stochastic increment
is bounded by It\^o's isometry and the linear-growth estimate:
\[
 \E\left|
 \int_s^t(\Pi_n\Sigma(q,X_n(q))^*\varphi_n,\dd W_q)_{\mathcal U}
 \right|
 \le C_\varphi |t-s|^{1/2}.
\]
This proves \eqref{2.9}.  It also yields the one-sided averaging bound
\[
 \left|\E[ZF_n(t)]-\frac1h\int_{t-h}^t\E[ZF_n(s)]\dd s\right|
 \le C_\varphi\|Z\|_\infty
      (h^{1/p}+h^{1/2}),\qquad 0<h<t,
\]
for every bounded random variable $Z$.  The weak convergence in
$L^2(\Omega\times(0,T);H_D)$ identifies the limit of the time average,
whereas the same estimate applies to the limiting variational equation.
Letting $h\downarrow0$ gives convergence against bounded $Z$.
Alternatively, and more directly, the fixed-time weak convergence proved
in Step~2 implies
\[
 \langle X_n(t),\varphi\rangle\rightharpoonup
 \langle X(t),\varphi\rangle\qquad\hbox{in }L^2(\Omega),
\]
because $\varphi\in H_0^1(\mathcal O)$ defines a continuous functional
on $H_D$.  This is precisely \eqref{2.8}. No strong convergence in
$L^2(\Omega;C([0,T];H_D))$ is used.
\end{proof}

We note that the regularity of time and space are not optimal, see \cite{BGW2022}.
\begin{remark}[Optional scalar smoothing]\label{rem2.9}
If a smooth scalar nonlinearity is desired, one may use
$\Phi_\varepsilon(r)=r(r^2+\varepsilon^2)^{(m-1)/2}+\varepsilon r$.
Then
\[
r\Phi_\varepsilon(r)\ge|r|^{m+1} \quad {\rm and} \quad
0<\Phi_\varepsilon'(r)\le C_m(1+|r|^{m-1}) \quad {\rm for} \quad 0<\varepsilon\le1.
\]
Its derivative is not globally bounded. The Galerkin construction and
energy estimates remain compatible, and monotone variational convergence
removes $\varepsilon$. The direct $C^1$ Galerkin construction above makes
this additional parameter unnecessary. In particular, the proof below does
not use an $\varepsilon$-regularization. If this optional route is adopted,
one must additionally pass $\varepsilon\downarrow0$ in the Malliavin
derivatives by a uniform derivative bound and the weak closedness of the
Malliavin derivative.
\end{remark}

\subsection{Malliavin differentiation before passage to the limit}

\begin{lemma}[Finite-dimensional derivative and its energy]
\label{lem2.10}
For fixed $n$, $X_n(t)$ is Malliavin differentiable. For $k\le n$ and
almost every $r\le t$, set $Y_{r,k}^n(s)=D_r^kX_n(s)$, $s\ge r$.
It satisfies
\begin{equation}\label{2.10}
 \begin{aligned}
 \dd Y_{r,k}^n&=P_n\Delta(a_nY_{r,k}^n)\dd s
      +\sum_{\ell\le n}P_nC_\ell Y_{r,k}^n \dd\beta_\ell(s),\\
 Y_{r,k}^n(r)&=P_n\Sigma(r,X_n(r))e_k,
 \quad a_n(s,x)=m|X_n(s,x)|^{m-1}.
 \end{aligned}
\end{equation}
The drift projection here is understood in the finite-dimensional
variational sense. With $\eta_{r,k}^n=Y_{r,k}^n(r)$,
\begin{equation}\label{2.11}
 \begin{aligned}
 \sup_{r\le s\le T}\E\|Y_{r,k}^n(s)\|_{H_D}^2
 +\E\int_r^T\!\int_{\mathcal O}a_n|Y_{r,k}^n|^2 \dd x\dd s
 \le C_T\E\|\eta_{r,k}^n\|_{H_D}^2.
 \end{aligned}
\end{equation}
Consequently, for $F_n=\langle X_n(t),\varphi\rangle$,
\begin{equation}\label{2.12}
 \sup_n\E\int_0^t\|D_rF_n\|_{\mathcal U}^2 \dd r<\infty.
\end{equation}
If $C=0$, the stronger pathwise inequality holds:
\begin{equation}\label{2.13}
 \frac12\|Y_{r,k}^n(s)\|_{H_D}^2
 +\int_r^s\!\int_{\mathcal O}a_n|Y_{r,k}^n|^2 \dd x\dd q
 \le\frac12\|P_ng(r)e_k\|_{H_D}^2.
\end{equation}
\end{lemma}

\begin{proof}
\emph{Step 1: verification of the finite-dimensional Malliavin theorem.}
Use $H_D$-orthonormal coordinates on $H_n$ and write
\[
 b_n(v)=-A_n(v),\qquad
 \sigma_n(s,v)=P_n\bigl(g(s)+Cv\bigr)\Pi_n.
\]
We verify the hypotheses of the deterministic-coefficient result
\cite[Corollary~3.5]{ImkellerReisSalkeld2019}, rather than obtaining
Malliavin differentiability by a state truncation. Since $m>1$, the scalar
map $\Phi(r)=|r|^{m-1}r$ is $C^1$. In fixed dimension this gives
$A_n\in C^1(H_n;H_n)$ and local Lipschitz continuity. Moreover,
\begin{align*}
 (b_n(v)-b_n(w),v-w)_{H_D}
 &=-\int_{\mathcal O}[\Phi(v)-\Phi(w)](v-w) \dd x\le0,\\
 (DA_n(v)y,y)_{H_D}
 &=\int_{\mathcal O}m|v|^{m-1}y^2 \dd x\ge0.
\end{align*}
Thus the drift is one-sided Lipschitz with constant zero. The diffusion is
$C^1$ in the state variable, with
\[
 D_v\sigma_n(s,v)y=P_nCy \Pi_n,
 \qquad
 \|\sigma_n(s,v)-\sigma_n(s,w)\|_{\mathcal L_2}^2
 \le L\|v-w\|_{H_D}^2
\]
by using \eqref{2.2}. Moreover, we have
\[
\int_0^T\|\sigma_n(s,0)\|_{\mathcal L_2}^2 \dd s<\infty
\]
because
$g$ is deterministic and continuous. The initial value $P_nx$ is
deterministic. Finally, the last assertion in Step~1 of
Proposition~\ref{prop2.8} gives, for every fixed $r\ge2$,
\[
 \sup_n\E\sup_{s\le T}(1+\|X_n(s)\|_{H_D}^2)^{r/2}<\infty.
\]
In particular the moment requirement with some exponent $r>2$ is satisfied.
Therefore all hypotheses of the cited finite-dimensional theorem hold and
$X_n$ is Malliavin differentiable. Its derivative is exactly \eqref{2.10}.
This argument is performed at fixed $n$, before any SPDE limit, and does not
require either a globally bounded $\Phi'$ or a state cutoff.

\smallskip
\emph{Step 2: relation to Cameron--Martin perturbations.}
For clarity, consider a deterministic control
$h\in L^2(0,T;\mathcal U)$ and the shift
\[
W^\theta_s=W_s+\theta\int_0^s h(r)\dd r.
\]
The shifted finite-dimensional equation acquires the extra drift
$\theta P_n\Sigma(s,X_n^\theta)\Pi_nh(s)$.
On a bounded set of controls, monotonicity, the Lipschitz noise estimate,
BDG, and Gronwall give
\[
 \E\sup_{s\le T}\|X_n^\theta(s)-X_n(s)\|_{H_D}^2
 \le C_{T,R}\theta^2\|h\|_{L^2(0,T;\mathcal U)}^2
       \bigl(1+\|x\|_{H_D}^2+\sup_{s\le T}\|g(s)\|_{\mathcal L_2}^2\bigr)
\]
when $\|h\|_{L^2}\le R$ and $|\theta|\le1$.
The drift in the difference quotient contains the averaged Jacobian
\[
\int_0^1 DA_n(X_n+\lambda(X_n^\theta-X_n)) \dd\lambda.
\]
In fixed dimension this converges locally to $DA_n(X_n)$, since
$A_n$ is $C^1$. Localization and the moment bounds identify the
directional derivative with the solution $Z_h$ of
\[
 \begin{split}
 \dd Z_h={}&-DA_n(X_n)Z_h \dd s
  +\sum_{\ell\le n}P_nC_\ell Z_h \dd\beta_\ell(s)
  +P_n\Sigma(s,X_n)\Pi_nh(s) \dd s,\qquad Z_h(0)=0.
 \end{split}
\]
Variation of constants, or stochastic Fubini followed by uniqueness of
this linear equation, then gives
\begin{equation}\label{2.14}
 D_hX_n(t)=Z_h(t)
       =\sum_{k\le n}\int_0^t Y_{r,k}^n(t)h_k(r) \dd r.
\end{equation}
In particular, an indicator control gives an integral of the kernel over
an interval. It does not give a value of $D_rX_n(t)$ at one time.
Equation \eqref{2.14} explains the derivative equation, while the
finite-dimensional theorem used in Step 1 supplies membership in the
Malliavin Sobolev space. Directional differentiability alone would not
supply that membership.

\smallskip
\emph{Step 3: estimates independent of dimension.}
Applying It\^o's formula to \eqref{2.10} gives
\[
 \begin{aligned}
 \dd\|Y\|_{H_D}^2+2\int_{\mathcal O}a_nY^2 \dd x\dd s
 &=\sum_{\ell\le n}\|P_nC_\ell Y\|_{H_D}^2\dd s
   +2\sum_{\ell\le n}(Y,P_nC_\ell Y)_{H_D}\dd\beta_\ell.
 \end{aligned}
\]
For simplicity, we write $Y=Y_{r,k}^n$ and $\eta=\eta_{r,k}^n$. We give the complete argument because the
uniformity in $n$, $r$ and $k$ is essential later.  For $R>0$, set
\[
 \tau_R=\inf\{s\ge r:\|Y(s)\|_{H_D}\ge R\}\wedge T.
\]
Integrating the preceding It\^o identity up to $s\wedge\tau_R$ and taking expectation eliminates the stopped local
martingale and gives, by \eqref{2.2},
\[
 \begin{aligned}
 &\E\|Y(s\wedge\tau_R)\|_{H_D}^2
 +2\E\int_r^{s\wedge\tau_R}\!\!\int_{\mathcal O}
       a_n|Y|^2\dd x\dd q\\
 &\le \E\|\eta\|_{H_D}^2
 +L\E\int_r^{s\wedge\tau_R}\|Y(q)\|_{H_D}^2\dd q .
 \end{aligned}
\]
Dropping the nonnegative space--time integral and applying Gronwall
yields
\[
 \sup_{r\le s\le T}
 \E\|Y(s\wedge\tau_R)\|_{H_D}^2
 \le e^{L(T-r)}\E\|\eta\|_{H_D}^2
 \le e^{LT}\E\|\eta\|_{H_D}^2.
\]
Returning to the preceding inequality with $s=T$ then gives
\[
 \begin{aligned}
 2\E\int_r^{T\wedge\tau_R}\!\!\int_{\mathcal O}
      a_n|Y|^2\dd x\dd q
 &\le \E\|\eta\|_{H_D}^2
 +L\int_r^T\E\|Y(q\wedge\tau_R)\|_{H_D}^2\dd q\\
 &\le (1+LT e^{LT})\E\|\eta\|_{H_D}^2.
 \end{aligned}
\]
The constants are independent of $n,r,k$ and $R$.  Letting
$R\uparrow\infty$, Fatou's lemma and monotone convergence give
\eqref{2.11}, with for instance
$C_T=e^{LT}+\frac12(1+LT e^{LT})$ after enlarging the constant.

\smallskip
The summation over the initial noise direction is legitimate because,
by orthogonality of $\Pi_n$ and contractivity of $P_n$,
\[
 \sum_{k\le n}\|\eta_{r,k}^n\|_{H_D}^2
 =\sum_{k\le n}
   \|P_n\Sigma(r,X_n(r))e_k\|_{H_D}^2
 \le\|\Sigma(r,X_n(r))\|_{\mathcal L_2(\mathcal U,H_D)}^2.
\]
Since $D_r^kF_n=\langle Y_{r,k}^n(t),\varphi\rangle$, duality and
\eqref{2.7} give the explicit bound
\[
 \E\int_0^t\|D_rF_n\|_{\mathcal U}^2 \dd r
 \le C_T\|\varphi\|_{H_0^1}^2
           \int_0^t\E\|\Sigma(r,X_n(r))\|_{\mathcal L_2}^2 \dd r.
\]
The right-hand side is uniform in $n$, proving \eqref{2.12}.
When $C=0$, both stochastic terms in the energy equation disappear;
integration gives \eqref{2.13} pathwise.
\end{proof}

\begin{lemma}[Weak closure at a fixed time]\label{lem2.11}
For every $t>0$ and $\varphi\in C_c^\infty(\mathcal O)$,
$F=\langle X(t),\varphi\rangle$ belongs to $\mathbb D^{1,2}$.
Along any weakly convergent subsequence of the derivatives in
\eqref{2.12}, the weak limit is $DF$.
\end{lemma}

\begin{proof}
The Malliavin derivative is a closed linear operator from $L^2(\Omega)$ to
$L^2(\Omega;L^2(0,t;\mathcal U))$; see \cite{Nualart2006}.
Its graph is a closed linear subspace of the product Hilbert space and
therefore is weakly closed. Now apply \eqref{2.8} and
\eqref{2.12}: along a subsequence,
$(F_n,DF_n)\rightharpoonup(F,V)$ in that product space. Weak closedness
implies $F\in\mathbb D^{1,2}$ and $V=DF$.

\smallskip
The identification can also be seen directly from Malliavin integration
by parts. For a smooth cylindrical random variable $Z$ and a deterministic
control $h$, write $\delta_{\!W}$ for the divergence operator. Then
\[
 \E[Z\langle DF_n,h\rangle_{L^2(0,t;\mathcal U)}]
       =\E[F_n\delta_{\!W}(Zh)].
\]
Both test quantities on the two sides are in $L^2$, so both pairings pass
under weak convergence. The resulting identity identifies $V$ as the
Sobolev derivative of $F$. Uniqueness of this derivative identifies every
weak subsequential limit. This step uses the scalar fixed-time convergence,
not an unproved strong path-space convergence of the approximate states.
\end{proof}

\subsection{The additive terminal estimate and full support}

Take $g(s)=G\in\mathcal L_2(\mathcal U,H_D)$ and $C=0$. The original
$L^2$-valued additive noise is included. Let $u=X$ and
$F_{t,\varphi}=\langle u(t),\varphi\rangle$.

\begin{proposition}[Additive terminal convergence]\label{prop2.12}
For $0<\delta\le\min(t,1)$,
\begin{equation}\label{2.15}
 \E\left[\frac1\delta\int_{t-\delta}^t
       \|D_rF_{t,\varphi}-G^*\varphi\|_{\mathcal U}^2 \dd r\right]
 \le C_{t,\varphi}\delta^{\frac{2}{m+1}}.
\end{equation}
\end{proposition}

\begin{proof}
\emph{Step 1: the finite-dimensional remainder.}
At the Galerkin level, pair \eqref{2.10} with $\varphi$ to obtain
\[
D_rF_n=\Pi_nG^*\varphi_n+R_r^n
\]
with
\[
R_{r,k}^n=\int_r^t\int_{\mathcal O}a_nY_{r,k}^n\Delta\varphi_n \dd x\dd s.
\]
Weighted Cauchy--Schwarz, followed by
\eqref{2.13} summed over $k$, yields
\[
 \|R_r^n\|_{\mathcal U}^2\le\tfrac12\|G\|_{\mathcal L_2}^2 C_{n,\delta},
 \qquad C_{n,\delta}=\int_{t-\delta}^t\!\int_{\mathcal O}
              a_n|\Delta\varphi_n|^2 \dd x\dd s.
\]
Indeed, the two factors before summing in $k$ are
\[
\int_r^t\int_{\mathcal O}a_n|Y_{r,k}^n|^2 \quad {\rm and} \quad
\int_r^t\int_{\mathcal O}a_n|\Delta\varphi_n|^2.
\]
The first has a pathwise bound from the additive derivative energy, and
the second is at most $C_{n,\delta}$ for $r\in[t-\delta,t]$.
The integral is over all of $\mathcal O$ because the projected test need
not retain the support of $\varphi$.

\smallskip
\emph{Step 2: smallness of the terminal strip.}
Set $\alpha=(m-1)/(m+1)$ and
\[
Z_n=\int_0^T\|X_n(s)\|_p^p \dd s.
\]
By \eqref{2.5},
\begin{equation}\label{2.16}
 C_{n,\delta}\le C_\varphi\delta^{1-\alpha}Z_n^\alpha,
 \qquad \sup_n\E Z_n<\infty.
\end{equation}
Thus
\[
\sup_n\E C_{n,\delta}\le C\delta^{1-\alpha}.
\]
To see the power of $\delta$, use H\"older on
$[t-\delta,t]\times\mathcal O$:
\[
 \int_{t-\delta}^t\!\int_{\mathcal O}|X_n|^{m-1} \dd x\dd s
 \le (\delta|\mathcal O|)^{1-\alpha}
       \left(\int_{t-\delta}^t\!\int_{\mathcal O}|X_n|^{m+1}
                                        \dd x\dd s\right)^\alpha.
\]
Jensen's inequality then uses only the first moment of $Z_n$. Notice that
$1-\alpha=2/(m+1)>0$ for every $m>1$; no nondegeneracy of the diffusion
coefficient $a_n$ is needed.

\smallskip
\emph{Step 3: passage to the degenerate equation.}
We finally justify the limit passage explicitly.  Fix $\delta>0$ and
work in the Hilbert space
\[
 \mathcal K_\delta
 =L^2\bigl(\Omega\times(t-\delta,t);\mathcal U\bigr).
\]
Lemma~\ref{lem2.11} gives, after extraction,
$DF_n\rightharpoonup DF$ in $\mathcal K_\delta$.  Since
$\varphi_n\to\varphi$ in $H_0^1$ and
$G\in\mathcal L_2(\mathcal U,H_D)$,
\[
 \Pi_nG^*\varphi_n\longrightarrow G^*\varphi
 \qquad\hbox{strongly \ in}\quad \mathcal U.
\]
Hence
\[
DF_n-\Pi_nG^*\varphi_n\rightharpoonup DF-G^*\varphi
\quad {\rm in} \quad \mathcal K_\delta.
\]
Weak lower semicontinuity of the squared
Hilbert norm therefore gives
\[
 \begin{aligned}
 \E\frac1\delta\int_{t-\delta}^t
   \|D_rF-G^*\varphi\|_{\mathcal U}^2\dd r
 &\le\liminf_{n\to\infty}
 \E\frac1\delta\int_{t-\delta}^t
   \|D_rF_n-\Pi_nG^*\varphi_n\|_{\mathcal U}^2\dd r\\
 &\le C_{t,\varphi}\delta^{1-\alpha}
 =C_{t,\varphi}\delta^{2/(m+1)}.
 \end{aligned}
\]
This is \eqref{2.15}.  The order of limits is important: first
$n\to\infty$ for fixed $\delta$, and only afterwards
$\delta\downarrow0$.
\end{proof}

For an adapted scalar observation at time $t$, define its Malliavin
covariance by
\[
\Gamma_F=\int_0^t\|D_rF\|_{\mathcal U}^2 \dd r.
\]
We need only positivity of this quantity, not an inverse-moment bound.
The following averaged criterion is useful because a Malliavin derivative
is an equivalence class in the time variable; pointwise convergence as
$r\uparrow t$ is unnecessary.

\begin{lemma}[Terminal covariance criterion]\label{lem2.13}
Suppose $F\in\mathbb D^{1,2}$ and $\Theta\in L^2(\Omega;\mathcal U)$ satisfy
\begin{equation}\label{2.17}
 \frac1\delta\int_{t-\delta}^t\|D_rF-\Theta\|_{\mathcal U}^2 \dd r
 \longrightarrow0\quad\hbox{in probability as }\delta\downarrow0.
\end{equation}
Then
\[
\{\|DF\|_{L^2(0,t;\mathcal U)}=0\}\subset\{\Theta=0\} \quad \hbox{up to
a null set.}
\]
In particular
\[
\Pp(F=0,\Theta\ne0)=0.
\]
If $\Theta\ne0$
almost surely, $F$ has an absolutely continuous law.
\end{lemma}

\begin{proof}
Select a deterministic sequence $\delta_j\downarrow0$ along which the
nonnegative averages in \eqref{2.17} converge almost
surely to zero. On $\{DF=0\}$ every such average equals $\|\Theta\|^2$.
The first inclusion follows. Locality of the Malliavin derivative gives
$DF=0$ almost surely on $\{F=0\}$, proving the second assertion.
Bouleau--Hirsch gives the density assertion when the covariance is positive
almost surely \cite{Nualart2006}.
\end{proof}

The locality step in Lemma~\ref{lem2.13} is stronger than
needed when $\Theta$ is deterministic and nonzero. It is retained because
the multiplicative terminal direction may vanish precisely on the
global zero-state event.
\begin{proposition}[Density of additive localized observations]
\label{prop2.14}
Under Assumption~\ref{ass2.4}, for every fixed $t>0$ and every
nonzero $\varphi\in C_c^\infty(\mathcal O)$, the law of
$F_{t,\varphi}$ is absolutely continuous with respect to Lebesgue measure.
In particular,
$\Pp(F_{t,\varphi}=a)=0$ for every $a\in\mathbb R$.
\end{proposition}

\begin{proof}
Lemma~\ref{lem2.11} gives the Sobolev derivative.
Proposition~\ref{prop2.12} implies
\eqref{2.17} with $\Theta=G^*\varphi$ by Markov's
inequality. This vector is nonzero by Assumption~\ref{ass2.4}.
The density conclusion follows from Lemma~\ref{lem2.13}.
A singleton is a Lebesgue-null set, which gives the last assertion.
\end{proof}

Absolute continuity in Proposition~\ref{prop2.14} does not assert
that the density is strictly positive everywhere, or that the law has
support all of $\mathbb R$. Only the absence of an atom at zero will be
used to rule out an interior spatial vacuum.

\begin{lemma}[Countable-basis support criterion]\label{lem2.15}
Let $v$ be a random distribution on $\mathcal O$. Fix a countable basis
$(O_j)$ of interior balls and an event $E$. Suppose that each $O_j$
admits a nonzero $\varphi_j\in C_c^\infty(O_j)$ such that
\[
\Pp(E\cap\{\langle v,\varphi_j\rangle=0\})=0.
\]
Then $\supp_{\mathcal O}v=\overline{\mathcal O}$ almost surely on $E$.
\end{lemma}

\begin{proof}
If the support is not full, the distribution vanishes on a nonempty
open subset of $\mathcal O$. That subset contains at least one basis
ball $O_j$, and hence $\langle v,\varphi_j\rangle=0$.
Consequently
\[
 E\cap\{\supp_{\mathcal O}v\ne\overline{\mathcal O}\}
 \subset\bigcup_{j\ge1}\bigl(E\cap\{\langle v,\varphi_j\rangle=0\}\bigr).
\]
The right-hand side is a countable union of null events. One can take
balls with rational centers and radii whose closures lie in $\mathcal O$;
this explains why an uncountable intersection over all open sets is not
required.
\end{proof}

\begin{theorem}[Additive full-domain support]\label{thm2.16}
Let $u_0\in H_D$ be deterministic and $G\in\mathcal L_2(\mathcal U,H_D)$.
Assume Assumption~\ref{ass2.4}. Then the solution of
\[
\dd u=\Delta\Phi(u)\dd t+G\dd W_t
\]
has, for every fixed $t>0$,
\begin{equation}\label{2.18}
 \Pp\bigl(\supp_{\mathcal O}u(t)=\overline{\mathcal O}\bigr)=1.
\end{equation}
Every nonzero smooth localized scalar observation has a density.
\end{theorem}

\begin{proof}
Proposition~\ref{prop2.14} gives the density assertion and excludes
a zero pairing for one chosen nonzero test in each basis ball.
Apply Lemma~\ref{lem2.15} with $E=\Omega$ to obtain
\eqref{2.18}.
\end{proof}

Proposition~\ref{prop2.5} supplies explicit trace-class
coefficients for Theorem~\ref{thm2.16}. The initial condition in this
theorem may be signed, or even zero. In particular, for the nonnegative
compactly supported data in Theorem~\ref{thm2.3},
\eqref{2.18} removes every interior vacuum at a fixed positive
time. This is instantaneous propagation in the precise sense of
Definition~\ref{def2.2}.

\section{Nonlocal multiplicative noise: the Dirichlet problem}

Let $H=H_D$ and let $B\in\mathcal L(H,\mathcal L_2(\mathcal U,H))$.
Write $B(v)e_k=K_kv$. We assume
\begin{equation}\label{3.1}
 \sum_{k\ge1}\|K_kv\|_H^2\le C_B\|v\|_H^2,\qquad v\in H,
\end{equation}
and consider
\begin{equation}\label{3.2}
 \dd u(t,x)=\Delta\Phi(u(t,x))\dd t+B(u(t,x))\dd W_t,t>0,\ x\in\mathcal{O}, \qquad u(0)=u_0\in H.
\end{equation}
This is the case $g=0$, $C=B$ of \eqref{2.3}; its solution
and Galerkin approximations are provided by Proposition~\ref{prop2.8}.
The deterministic initial condition may have either sign.

\subsection{Well-posedness and nonlocal vacuum ellipticity}
Similar to the additive noise, we first introduce the following definition.
\begin{definition}\label{def3.1}
A family $(K_k)$ is vacuum preserving if, for every open
$O\subset\mathcal O$, $v|_O=0$ implies $(K_kv)|_O=0$ for all $k$.
The restrictions are understood in the distributional sense.
\end{definition}

Pointwise multiplication by sufficiently regular functions has the
property in Definition~\ref{def3.1}. An integral operator need not:
values of $v$ outside $O$ can contribute to $K_kv$ inside $O$.
This is the distinction relevant to activation of a spatial vacuum.
It is separate from well-posedness, for which the Hilbert--Schmidt
bound \eqref{3.1} is sufficient.

\begin{proposition}[Well-posedness]\label{prop3.2}
Under \eqref{3.1}, equation \eqref{3.2} has a unique
variational solution $u$ in
$L^2(\Omega;C([0,T];H_D))\cap L^p(\Omega\times(0,T)\times\mathcal O)$.
For each finite $T$,
\[
 \E\sup_{s\le T}\|u(s)\|_{H_D}^2
       +\E\int_0^T\|u(s)\|_p^p \dd s
 \le C_T\|u_0\|_{H_D}^2.
\]
The zero initial condition gives the identically zero solution.
\end{proposition}

\begin{proof}
Apply Proposition~\ref{prop2.8} with $g=0$ and $C=B$.
For completeness, the energy balance now has no additive forcing:
\[
 \dd\|u\|_{H_D}^2+2\|u\|_p^p \dd s
  =\|B(u)\|_{\mathcal L_2}^2 \dd s
       +2(u,B(u)\dd W_s)_{H_D}.
\]
The same BDG and Gronwall argument therefore gives the displayed bound
with right-hand side proportional to $\|u_0\|_{H_D}^2$.
For two solutions, monotonicity and \eqref{3.1} control the squared
$H_D$ norm of their difference, proving uniqueness. Since both drift and
noise vanish at zero, uniqueness also identifies the zero solution.
\end{proof}

The additional assumption below addresses which interior tests see the
noise at a nonzero state.

\begin{assumption}[Nonlocal vacuum ellipticity]\label{ass3.3}
For every $v\in H\setminus\{0\}$ and every nonzero
$\varphi\in C_c^\infty(\mathcal O)$,
\begin{equation}\label{3.3}
 \|B(v)^*\varphi\|_{\mathcal U}^2
 =\sum_{k\ge1}|\langle K_kv,\varphi\rangle|^2>0.
\end{equation}
\end{assumption}

Assumption~\ref{ass3.3}
forces every nonzero state to be visible through every nonzero interior
test. In particular, if $v$ vanishes on $O$ but is globally nonzero,
no nonzero test supported in $O$ is orthogonal to all the noise directions.
Thus nonlocal vacuum ellipticity is incompatible with vacuum preservation for such states.
It is a qualitative strict-positivity assumption, with no uniform
ellipticity constant as $v$ or the test varies.

\begin{proposition}[An explicit coefficient]
\label{prop3.4}
Let $(h_j)$ be an orthonormal basis of $H$. On the noise space indexed by
$(i,j)\in\mathbb N^2$, define
\begin{equation*}
 K_{ij}v=2^{-i-j}(v,h_j)_Hh_i.
\end{equation*}
Then \eqref{3.1} and \eqref{3.3} hold.
\end{proposition}

\begin{proof}
Since $\sum_{i\ge1}4^{-i}=1/3$,
\[
\sum_{i,j}\|K_{ij}v\|_H^2
\le\frac{1}{12}\|v\|_H^2.
\]
More explicitly,
\[
 \sum_{i,j}\|K_{ij}v\|_H^2
 =\left(\sum_i4^{-i}\right)
         \sum_j4^{-j}|(v,h_j)_H|^2
 \le\frac1{12}\sum_j|(v,h_j)_H|^2.
\]
If $v\ne0$, some $(v,h_j)_H$ is nonzero.
If $\varphi\ne0$, its continuous functional on $H$ is nonzero, so
some $\langle h_i,\varphi\rangle$ is nonzero. The corresponding term in
\eqref{3.3} is positive.
For the Dirichlet eigensystem one may take $h_j=\sqrt{\lambda_j}\rho_j$;
then for function-valued $v$,
\[
K_{ij}v=2^{-i-j}\sqrt{\lambda_i/\lambda_j}
\rho_i\int_{\mathcal O}v\rho_j \dd x,
\]
an explicit rank-one integral operator.
\end{proof}

The two indices in Proposition~\ref{prop3.4} have different roles:
$j$ detects the nonzero state and $i$ distributes a noise direction
visible to the chosen spatial test. One can more generally use
$K_{ij}v=a_{ij}\ell_j(v)f_i$, where the bounded functionals $\ell_j$
separate points, the profiles $f_i$ separate nonzero interior tests,
$a_{ij}>0$, and
\[
\sum_{i,j}a_{ij}^2\|\ell_j\|^2\|f_i\|_H^2<\infty.
\]
The same calculation proves the Hilbert--Schmidt bound and nonlocal vacuum ellipticity.
Neither construction by itself proves that a solution cannot reach zero.

\subsection{Malliavin differentiation for the multiplicative equation}

It is useful to make explicit which parts of the additive construction
survive. Let $u_n$ solve \eqref{2.6} with $g=0$ and $C=B$.
For $r\le s$, its derivative has initial value $P_nK_ku_n(r)$ and satisfies
\[
 \dd Y_{r,k}^n=P_n\Delta(m|u_n|^{m-1}Y_{r,k}^n) \dd s
                 +\sum_{\ell\le n}P_nK_\ell Y_{r,k}^n \dd\beta_\ell(s).
\]
In particular, the stochastic term remains in the linearized equation.

\begin{proposition}[Multiplicative scalar Malliavin bound]
\label{prop3.5}
Under \eqref{3.1}, for each fixed $t>0$ and
$\varphi\in C_c^\infty(\mathcal O)$,
$F_{t,\varphi}=\langle u(t),\varphi\rangle$ belongs to $\mathbb D^{1,2}$.
The approximate derivatives satisfy
\[
 \sup_n\E\int_0^t\|D_r\langle u_n(t),\varphi\rangle\|_{\mathcal U}^2 \dd r
 \le C_{T,\varphi}C_B\int_0^t
                      \sup_n\E\|u_n(r)\|_{H_D}^2 \dd r<\infty.
\]
\end{proposition}

\begin{proof}
Lemma~\ref{lem2.10} applies to these finite-dimensional
coefficients. Summing \eqref{2.11} over the initial noise
directions gives
\[
 \begin{split}
 &\sup_{r\le s\le T}\sum_{k\le n}\E\|Y_{r,k}^n(s)\|_{H_D}^2
  +\sum_{k\le n}\E\int_r^T\!\int_{\mathcal O}
                       m|u_n|^{m-1}|Y_{r,k}^n|^2 \dd x\dd s\\
 &\le C_T\E\|B(u_n(r))\|_{\mathcal L_2}^2
                     \le C_TC_B\E\|u_n(r)\|_{H_D}^2.
 \end{split}
\]
For completeness, the scalar estimate follows directly from duality.
Since
$D_r^k\langle u_n(t),\varphi\rangle
 =\langle Y_{r,k}^n(t),\varphi\rangle$,
\[
 \begin{aligned}
 \E\|D_r\langle u_n(t),\varphi\rangle\|_{\mathcal U}^2=\sum_{k\le n}\E
      |\langle Y_{r,k}^n(t),\varphi\rangle|^2
 &\le \|\varphi\|_{H_0^1}^2
      \sum_{k\le n}\E\|Y_{r,k}^n(t)\|_{H_D}^2\\
 &\le C_T\|\varphi\|_{H_0^1}^2
      \E\|B(u_n(r))\|_{\mathcal L_2}^2\\
 &\le C_T C_B\|\varphi\|_{H_0^1}^2
      \E\|u_n(r)\|_{H_D}^2.
 \end{aligned}
\]
Integrating in $r\in(0,t)$ and using the uniform state estimate
\eqref{2.7} gives
\[
 \sup_n\E\int_0^t\|D_r\langle u_n(t),\varphi\rangle\|_{\mathcal U}^2\dd r
 \le C_{T,\varphi}C_B\int_0^t
                      \sup_n\E\|u_n(r)\|_{H_D}^2\dd r<\infty.
\]
The scalar weak convergence in Proposition~\ref{prop2.8} and
the graph argument in Lemma~\ref{lem2.11} then identify the
derivative of the limiting observable.
\end{proof}

Proposition~\ref{prop3.5} provides the Sobolev regularity
needed below. Its derivative energy is an expectation estimate.
It cannot be inserted as a deterministic bound into a pathwise product
with the random coefficient $m|u_n|^{m-1}$.

Next, we consider the random remainder in a mixed norm. Note that the pathwise contraction \eqref{2.13} is unavailable
when the derivative equation contains stochastic terms. We instead use
\eqref{2.11} in expectation and a mixed norm. Let
\begin{equation}\label{3.4}
 \alpha=\frac{m-1}{m+1},\qquad
 q=\frac{m+1}{2m}=\frac1{1+\alpha},\qquad
 \rho=2q=\frac{m+1}{m}\in(1,2).
\end{equation}
In particular, $q(1-\alpha)=1/m$.

\begin{proposition}[Affine-noise terminal convergence]
\label{prop3.6}
For the affine equation \eqref{2.3}, fix $t>0$ and
$\varphi\in C_c^\infty(\mathcal O)$, and set
$F=\langle X(t),\varphi\rangle$ and
$\Theta=\Sigma(t,X(t))^*\varphi$. Then $F\in\mathbb D^{1,2}$ and
\begin{equation}\label{3.5}
 \E\left[\left(\frac1\delta\int_{t-\delta}^t
       \|D_rF-\Theta\|_{\mathcal U}^2 \dd r\right)^q\right]
 \longrightarrow0,\qquad  as \ \delta\downarrow0.
\end{equation}
Only the second-order state and linearized estimates are required.
\end{proposition}

\begin{proof}
We divide the proof into five steps.

\smallskip
\emph{Step 1: retain the approximate state at the differentiation time.}
At the Galerkin level write
\begin{equation}\label{3.6}
 D_rF_n-b_n(r)=R_r^n+S_r^n,
 \qquad b_n(r)=\Pi_n\Sigma(r,X_n(r))^*\varphi_n,
\end{equation}
where, for $k\le n$,
\[
 \begin{aligned}
 R_{r,k}^n=\int_r^t\!\int_{\mathcal O}
                  a_nY_{r,k}^n\Delta\varphi_n \dd x\dd s,\quad
 S_{r,k}^n=\sum_{\ell\le n}\int_r^t
                  \langle C_\ell Y_{r,k}^n(s),\varphi_n\rangle
                   \dd\beta_\ell(s).
 \end{aligned}
\]
The center in \eqref{3.6} is evaluated at $r$, not at $t$.
This choice removes any need for a uniform time modulus of the approximate
states in $H$.

\smallskip
\emph{Step 2: control the product in the drift remainder.}
Define
\[
A_n(r)=\sum_{k\le n}\int_r^t\int_{\mathcal O}a_n|Y_{r,k}^n|^2 \dd x\dd s\quad
{\rm and} \quad \overline A_{n,\delta}=\delta^{-1}\int_{t-\delta}^t A_n(r) \dd r.
\]
Let $C_{n,\delta}$ and $Z_n$ be as in
\eqref{2.16}. Weighted Cauchy--Schwarz gives
\[
 \frac1\delta\int_{t-\delta}^t\|R_r^n\|_{\mathcal U}^2 \dd r
 \le\overline A_{n,\delta} C_{n,\delta}.
\]
Equations \eqref{2.11} and \eqref{2.7} yield
$\sup_{n,\delta}\E\overline A_{n,\delta}<\infty$.
By \eqref{2.16},
\[
\E C_{n,\delta}^{1/\alpha}
\le C\delta^{(1-\alpha)/\alpha}\E Z_n.
\]
H\"older with exponents $1/q$ and $1/(1-q)$ therefore gives the valid
random-product estimate
\begin{equation}\label{3.7}
 \begin{aligned}
 \E(\overline A_{n,\delta}C_{n,\delta})^q
\le(\E\overline A_{n,\delta})^q
       (\E C_{n,\delta}^{1/\alpha})^{\alpha q}\le C\delta^{q(1-\alpha)}=C\delta^{\frac1m}.
 \end{aligned}
\end{equation}
In particular, no expectation bound has been substituted into a pathwise
inequality, and no product of dependent random variables has been factored.
The choice of $q$ is exactly what makes this possible:
$q/(1-q)=1/\alpha$. Thus the required moment of $C_{n,\delta}$ reduces
to the first moment of $Z_n$, already available from the state energy.
An estimate of $\E(\overline A_{n,\delta}C_{n,\delta})$ would require
additional integrability that has not been proved.

\smallskip
\emph{Step 3: control the stochastic remainder.}
For the stochastic remainder, It\^o isometry, duality and
\eqref{2.2} give
\[
 \begin{aligned}
 \E\|S_r^n\|_{\mathcal U}^2
 &\le\|\varphi_n\|_{H_0^1}^2
       \sum_{k\le n}\int_r^t\E\sum_{\ell\le n}
                       \|C_\ell Y_{r,k}^n(s)\|_H^2 \dd s\\
 &\le C\int_r^t\E\|\Sigma(r,X_n(r))\|_{\mathcal L_2}^2 \dd s
 \le C(t-r).
 \end{aligned}
\]
Jensen's inequality and \eqref{3.7} imply
\begin{equation}\label{3.8}
 \sup_n\E\left[\left(\frac1\delta\int_{t-\delta}^t
          \|D_rF_n-b_n(r)\|_{\mathcal U}^2 \dd r\right)^q\right]
 \le C(\delta^{\frac1m}+\delta^q).
\end{equation}

Here $0<q<1$, so the passage from the separate remainders uses
$(a+b)^q\le a^q+b^q$ for nonnegative $a,b$, after the squared triangle
inequality. For the stochastic term,
\[
\E[(\delta^{-1}\int\|S_r^n\|^2 \dd r)^q]
\le(\delta^{-1}\int\E\|S_r^n\|^2 \dd r)^q\le C\delta^q.
\]

\emph{Step 4: pass to the weak limit in the appropriate Banach space.}
Fix $\delta$ before passing to the limit. Weak convergence of $X_n$ in
$L^2(\Omega\times(0,T);H)$ and affine linearity of $\Sigma$ imply
$b_n\rightharpoonup b$, where $b(r)=\Sigma(r,X(r))^*\varphi$.
The moving projections cause no difficulty: they are uniformly bounded
and converge strongly when applied to each fixed dual test.
By Lemma~\ref{lem2.11}, $DF_n\rightharpoonup DF$ in $L^2$.
Because $\rho=2q>1$, the space
$L^\rho(\Omega;L^2(t-\delta,t;\mathcal U))$ is a Banach space and its
norm is weakly lower semicontinuous. The continuous embedding of $L^2$
into this space permits passage to the limit in
\eqref{3.8}. Thus the same estimate holds for $DF-b$.
This step uses a norm with exponent $\rho>1$, not convexity of $x\mapsto x^q$.
Indeed the left-hand side of \eqref{3.8} is
\[
\delta^{-q}\|DF_n-b_n\|_{L^\rho(\Omega;L^2(t-\delta,t;\mathcal U))}^{\rho}.
\]
The inclusion from the ambient $L^2$ space is continuous because the
probability measure has mass one and $\rho<2$. It is therefore weakly
continuous, and the usual norm lower semicontinuity applies.

\smallskip
\emph{Step 5: move the center only for the limiting solution.}
Finally, continuity of $g$ and
$X\in L^2(\Omega;C([0,T];H))$ give
\[
 \omega_t(\delta):=\frac1\delta\int_{t-\delta}^t
       \E\|b(r)-b(t)\|_{\mathcal U}^2 \dd r\longrightarrow0.
\]
The averaged error with center $\Theta=b(t)$ has $q$th moment at most
\[
C(\delta^{1/m}+\delta^q+\omega_t(\delta)^q).
\]
This proves \eqref{3.5}. The order of operations is:
pass $n\to\infty$ with center $b_n(r)$, replace $b(r)$ by $b(t)$ for
the limiting solution, and then let $\delta\downarrow0$.
\end{proof}

\subsection{Full support or the zero state}

The terminal estimate has now been proved for the coefficient class
\eqref{3.1}. The remaining issue is dynamical: multiplicative noise
vanishes at the global zero state. Nonlocal vacuum ellipticity distinguishes every nonzero state
from that state but does not itself establish a non-hitting theorem.

\begin{theorem}\label{thm3.7}
Assume \eqref{3.1} and Assumption~\ref{ass3.3}, and let $u_0\in H_D$
be deterministic. For each fixed $t>0$,
\begin{equation}\label{3.9}
 \Pp\bigl(u(t)=0\ \hbox{or}\
          \supp_{\mathcal O}u(t)=\overline{\mathcal O}\bigr)=1.
\end{equation}
The two events are disjoint. In particular,
\begin{equation}\label{3.10}
 \Pp\bigl(\supp_{\mathcal O}u(t)=\overline{\mathcal O}\bigr)
 =\Pp\bigl(u(t)\ne0\bigr).
\end{equation}
For every nonzero $\varphi\in C_c^\infty(\mathcal O)$,
\[
\Pp(\langle u(t),\varphi\rangle=0,\ u(t)\ne0)=0.
\]
\end{theorem}

\begin{proof}
Apply Proposition~\ref{prop3.6} with $g=0$ and $C=B$.
Markov's inequality yields \eqref{2.17} with
$\Theta=B(u(t))^*\varphi$. By Assumption~\ref{ass3.3}, $\Theta=0$ implies $u(t)=0$.
Lemma~\ref{lem2.13} and locality therefore give
\[
\Pp(F_{t,\varphi}=0,u(t)\ne0)=0.
\]
Apply Lemma~\ref{lem2.15} with $E=\{u(t)\ne0\}$.
This retains the same survival event for all the countably many tests and
proves \eqref{3.9}.
Since the zero distribution has empty support, the disjointness and
\eqref{3.10} follow.
\end{proof}
If we want to get the full-support, we need the following definition.
\begin{definition}[Backward non-hitting for an initial state]\label{def3.8}
For a deterministic $u_0\ne0$, the solution has the backward non-hitting
property, abbreviated BNH, if
\begin{equation}\label{3.11}
 \Pp(\tau_0=\infty)=1,\qquad
 \tau_0=\inf\{s\ge0:u(s)=0\text{ in }H_D\}.
\end{equation}
The infimum of the empty set is $\infty$. A coefficient family has BNH
for a class of initial data if this property holds for every nonzero
initial datum in that class.
\end{definition}

Continuity in $H_D$ makes the hitting time in Definition~\ref{def3.8}
well defined. Since zero is absorbing, BNH rules out finite-time
extinction from the chosen nonzero initial state. It is stronger than
the nonvanishing assertion needed at one particular deterministic time.

\begin{corollary}[Dirichlet full support under non-hitting]
\label{cor3.9}
Under the hypotheses of Theorem~\ref{thm3.7}, suppose in addition
that $u_0\ne0$ and BNH holds in the sense of Definition~\ref{def3.8}.
Then $u(t)$ has full-domain support almost surely at each fixed $t>0$,
and every nonzero localized scalar observation has a density.
For a conclusion at one prescribed time, the weaker assumption
$\Pp(u(t)\ne0)=1$ suffices.
\end{corollary}

\begin{proof}
Use \eqref{3.10}. For the density assertion, nonlocal vacuum ellipticity makes the
terminal direction nonzero almost surely; apply the density part of
Lemma~\ref{lem2.13}.
\end{proof}

\begin{remark}[The remaining Dirichlet question]\label{rem3.10}
The present estimates do not establish \eqref{3.11} for every family
satisfying \eqref{3.1} and nonlocal vacuum ellipticity, including the example in
Proposition~\ref{prop3.4}. To see the issue in a logarithmic
argument, put $R=\|u\|_H^2$ and stop before $R$ reaches zero. It\^o's
formula gives
\[
 \begin{split}
 \dd\log R={}&\left(
 -\frac{2\|u\|_p^p}{R}
 +\frac{\sum_k\|K_ku\|_H^2}{R}
 -\frac{2\sum_k(u,K_ku)_H^2}{R^2}\right)\dd t+\frac2R\sum_k(u,K_ku)_H \dd\beta_k.
 \end{split}
\]
The noise drift and the quadratic variation rate of the martingale are
bounded by \eqref{3.1}. The negative drift
$-2\|u\|_p^p/R$, however, is not bounded from below by that condition.
The energy estimate controls $\int\|u\|_p^p$ but not
$\int\|u\|_p^p/\|u\|_H^2$. It therefore does not exclude
$\log R\to-\infty$ in finite time.
This is a missing non-hitting argument, not a counterexample to non-hitting.
Corollary~\ref{cor3.9} explicitly keeps that extra input; the unconditional
conclusion proved here is Theorem~\ref{thm3.7}.
\end{remark}

An invertible flow for the linear noise equation would not alone settle
Remark~\ref{rem3.10} either. Such a transformation removes the
stochastic differential, but one still needs injectivity or an appropriate
backward-uniqueness property for the transformed nonlinear evolution.
Forward uniqueness and the absorbing nature of zero do not imply that
a nonzero trajectory cannot arrive there.

\smallskip
The mass mechanism below applies to a different boundary problem. The
constant function $1$ does not belong to $H_0^1(\mathcal O)$, so mass is
not a continuous functional on $H_D$ and cannot simply be imported into
the Dirichlet argument in Remark~\ref{rem3.10}.

\section{Nonlocal multiplicative noise: the Neumann problem}

We now construct a conservative model for which the mass explicitly
excludes the global zero state. Its Malliavin argument uses only the
zero-mean noise coordinates. This avoids imposing a terminal estimate as
an additional hypothesis.

\smallskip
The equation has the boundary interpretation
\begin{equation}\label{4.1}
 \begin{cases}
 \dd u(t,x)=\Delta\Phi(u(t,x)) \dd t+B_N(u(t,x)) \dd W_t,
       &(t,x)\in(0,T)\times\mathcal O,\\
 \partial_\nu\Phi(u(t,x))=0,
       &(t,x)\in(0,T)\times\partial\mathcal O,\\
 u(0,x)=u(x)_0, & x\in\mathcal O.
 \end{cases}
\end{equation}
The Neumann boundary condition removes the diffusion flux through
$\partial\mathcal O$. We specify below a coefficient $B_N$ whose mass
coordinate can be solved explicitly. In contrast with the Dirichlet
state space, integration against the constant function is continuous in
the Neumann state space we use.

\subsection{Zero-mean equation and well-posedness}
Let
\[
V_N^1=\{\psi\in H^1(\mathcal O):\int_{\mathcal O}\psi \dd x=0\},
\]
equipped with the gradient norm, and let $H_0=(V_N^1)^*$.
An element of $H_0$ acts on $H^1$ by first subtracting the mean, so its
mass is zero. With $A_{N,0}=-\Delta_N$ on positive Neumann modes, set
\[
\|z\|_{H_0}^2=\langle z,A_{N,0}^{-1}z\rangle.
\]
Use the full state space
\begin{equation*}
 \mathcal H_N=\mathbb R\oplus H_0,\quad
 v=M(v)f_0+v^0,\quad f_0=|\mathcal O|^{-1}{\bf1}_{\mathcal O},\quad
 \|v\|_{\mathcal H_N}^2=|M(v)|^2+\|v^0\|_{H_0}^2.
\end{equation*}
Thus $M$ is the mass coordinate, equal to $\int v$ when $v$ is a function.
The boundary condition is $\partial_\nu\Phi(u)=0$ in the Neumann
variational sense.

\begin{proposition}[Monotonicity on a fixed-mass space]
\label{prop4.1}
If $u=Mf_0+z$ and $v=Mf_0+w$ have the same mass, then
\[
 \langle A_M(z)-A_M(w),z-w\rangle
 =\int_{\mathcal O}[\Phi(u)-\Phi(v)](u-v) \dd x\ge0,
\]
where
\[
\langle A_M(z),\psi\rangle
=\int_{\mathcal O}\Phi(z+M/|\mathcal O|)\psi \dd x
\]
for zero-mean $\psi$. For two states with unequal masses the analogous
expression need not be nonnegative.
\end{proposition}

\begin{proof}
The difference $z-w=u-v$ has zero mass. If $\xi$ is any zero-mean
distribution, its squared $H_0$ norm is
$\langle\xi,A_{N,0}^{-1}\xi\rangle$. Hence the drift paired with
$z-w$ in this norm becomes the displayed integral. Since $\Phi$ is
nondecreasing, the integrand is nonnegative pointwise for equal mass.
For unequal masses the zero-mean projection of $u-v$ replaces $u-v$
in the integral. It is not pointwise aligned with $\Phi(u)-\Phi(v)$,
as the calculation below shows.
\end{proof}

For example, take $m=2$ on a domain of unit volume, split into two
halves, with the constant values $u=(10,2)$ and $v=(9,0)$ on those
halves. Then $u-v=(1,2)$ has mean $3/2$ and zero-mean part
$(-1/2,1/2)$. Since $\Phi(r)=r|r|$, the putative monotonicity integral is
\[
 \frac12(100-81)\times\Big(-\frac12\Big)+\frac12(4-0)\times\frac12
 =-\frac{15}{4}.
\]
This explains why treating the full mass-plus-zero-mean product as a
single monotone $H^{-1}$ equation would be invalid. We solve the mass
first and apply Proposition~\ref{prop4.1} only within
its resulting fixed-mass equation.

\smallskip
Choose $f_i\in L^2_0(\mathcal O)$ whose span is dense in $L^2_0$, and
bounded functionals $\ell_j\in H_0^*$ that separate points of $H_0$.
For $\sigma\in\mathbb R$, $b_i>0$ and $a_{ij}>0$, define
\begin{equation}\label{4.2}
 K_0v=\sigma M(v)f_0,\qquad
 K_{i0}v=b_iM(v)f_i,\qquad
 K_{ij}v=a_{ij}\ell_j(v^0)f_i\quad(i,j\ge1).
\end{equation}
Assume
\begin{equation}\label{4.3}
 \sum_i b_i^2\|f_i\|_{H_0}^2<\infty,
 \qquad \sum_{i,j}a_{ij}^2\|\ell_j\|_{H_0^*}^2\|f_i\|_{H_0}^2<\infty.
\end{equation}
All driving Brownian motions $\beta_0$, $(\beta_{i0})$ and
$(\beta_{ij})$ are independent. Let $W'$ denote the cylindrical noise
formed by the last two families. The equation is
\begin{equation}\label{4.4}
 \dd u=\Delta_N\Phi(u)\dd t+K_0u \dd\beta_0
       +\sum_i K_{i0}u \dd\beta_{i0}
       +\sum_{i,j}K_{ij}u \dd\beta_{ij}.
\end{equation}
The sums in \eqref{4.3} make the coefficient a bounded linear map
from $\mathcal H_N$ to the relevant Hilbert--Schmidt space.
Indeed,
\[
 \|B_N(v)\|_{\mathcal L_2}^2
 =\sigma^2M(v)^2+M(v)^2\sum_i b_i^2\|f_i\|_{H_0}^2
        +\sum_{i,j}a_{ij}^2|\ell_j(v^0)|^2\|f_i\|_{H_0}^2,
\]
which is at most a constant times $\|v\|_{\mathcal H_N}^2$.
All modes with indices $i0$ or $ij$ have zero mass; only the mode
indexed by $0$ changes the mass.

\begin{proposition}[A concrete Neumann family]\label{prop4.2}
The family \eqref{4.2} satisfies nonlocal vacuum ellipticity on $\mathcal H_N$ against
every nonzero $\varphi\in C_c^\infty(\mathcal O)$. An explicit choice is
\[
f_i=h_i, \quad \ell_j(z)=(z,h_j)_{H_0}, \quad b_i=2^{-i} \quad {\rm and} \quad a_{ij}=2^{-i-j},
\]
where $h_j=\sqrt{\lambda_j}\rho_j$ for the positive $L^2$-normalized
Neumann eigenfunctions.
\end{proposition}

\begin{proof}
A nonzero compactly supported smooth test is not constant on $\mathcal O$.
Density of $(f_i)$ therefore supplies an $i$ with
$\langle f_i,\varphi\rangle\ne0$. If $M(v)\ne0$, then
\[
\langle K_{i0}v,\varphi\rangle\ne0.
\]
If $M(v)=0$ and $v\ne0$, choose
$j$ with $\ell_j(v^0)\ne0$, then
\[
\langle K_{ij}v,\varphi\rangle\ne0.
\]
This proves nonlocal vacuum ellipticity.
For the explicit choice, $(h_i)$ is orthonormal in $H_0$ and its span is
dense in $L^2_0$. Each functional has norm one, and the two sums in
\eqref{4.3} are $1/3$ and $1/9$, respectively.
\end{proof}

For the nonzero-mass class, the $K_{i0}$ modes alone already make every
nonzero interior test visible: their profiles are dense among zero-mean
functions, and a nonzero compactly supported test cannot be constant
on the connected domain. The $K_{ij}$ modes extend nonlocal vacuum ellipticity to nonzero
zero-mass states. This distinction will matter when we describe the
scope of the support theorem.

\smallskip
Define
\[G_0:\ell^2(\mathbb N)\to H_0, \quad G_0e_i=b_if_i
\quad
{\rm and}
\quad C:H_0\to\mathcal L_2(\ell^2(\mathbb N^2),H_0), \quad
(Cz)e_{ij}=a_{ij}\ell_j(z)f_i.
\]
In formulas involving $W'$, extend
these operators by zero to the other noise coordinates.

\begin{proposition}[Well-posedness]
\label{prop4.3}
For every deterministic $u_0=M_0f_0+z_0\in\mathcal H_N$,
equation \eqref{4.4} has a unique global variational solution
with continuous $\mathcal H_N$ paths. Its mass is
\begin{equation}\label{4.5}
 M_t=M_0\exp\{\sigma\beta_0(t)-\tfrac12\sigma^2t\},
 \qquad u(t)=M_tf_0+z(t).
\end{equation}
For every finite $T$ the solution satisfies
\[
\E\sup_{s\le T}\|u(s)\|_{\mathcal H_N}^2+
\E\int_0^T\|u(s)\|_p^p \dd s<\infty.
\]
If $M_0\ne0$, then $u(t)\ne0$ for every finite $t$, on one event of
probability one.
\end{proposition}

\begin{proof}
Testing the equation by the constant $1$ gives
\[
\dd M_t=\sigma M_t\dd\beta_0(t),
\]
since every other noise mode and
the Neumann drift have zero mass. Its unique solution is
\eqref{4.5}. Conversely, start with this mass process and put
$c(t)=M_t/|\mathcal O|$. Construct $z$ in the fixed triple
\[
L^p_0(\mathcal O)\subset H_0\subset (L^p_0(\mathcal O))^*
\]
from
\begin{equation}\label{4.6}
 \dd z=\Delta_N\Phi(z+c(t))\dd t
             +M_tG_0\dd W'_t+Cz\dd W'_t,\qquad z(0)=z_0.
\end{equation}
Here $c(t)$ is a spatial constant and the two displayed stochastic terms
act in disjoint coordinates of $W'$.
For any fixed real $c$, define $A_c$ on $L^p_0$ by
\[
\langle A_c(z),w\rangle=\int_{\mathcal O}\Phi(z+c)w \dd x.
\]
For $z,w\in L^p_0$,
\begin{equation}\label{4.7}
 \begin{aligned}
 &\langle A_c(z)-A_c(w),z-w\rangle
 =\int_{\mathcal O}[\Phi(z+c)-\Phi(w+c)](z-w) \dd x\ge0,\\
 &\qquad \langle A_c(z),z\rangle
 =\int_{\mathcal O}[\Phi(z+c)-\Phi(c)]z \dd x
 \ge2^{1-m}\|z\|_p^p.
 \end{aligned}
\end{equation}
The last inequality is the scalar uniform monotonicity of
$r\mapsto|r|^{p-2}r$. Also
\[
\|A_c(z)\|_{(L^p_0)^*}^{p/(p-1)}\le C(\|z\|_p^p+|c|^p).
\]
Since $p=m+1\ge2$, the standard
uniform monotonicity inequality
\[
 \bigl(|a|^{p-2}a-|b|^{p-2}b\bigr)(a-b)
 \ge 2^{2-p}|a-b|^p,\qquad a,b\in\mathbb R,
\]
applied with $a=z(x)+c$ and $b=c$, gives pointwise
\[
 [\Phi(z+c)-\Phi(c)]z
 \ge2^{2-p}|z|^p=2^{1-m}|z|^p.
\]
Since $z$ has zero mean,
$\int_{\mathcal O}\Phi(c)z\dd x=0$, and integration proves the
second line of \eqref{4.7}.  The estimate is valid for every real
$c$ and every signed $z$.

\smallskip
The growth estimate follows from the embedding of the zero-mean dual
into the full $L^{p'}$ dual:
\[
 \begin{aligned}
 \|A_c(z)\|_{(L^p_0)^*}
 &\le\|\Phi(z+c)\|_{p'}
   =\|z+c\|_p^{p-1},\\
 \|A_c(z)\|_{(L^p_0)^*}^{p'}
 &\le\|z+c\|_p^p
 \le2^{p-1}\bigl(\|z\|_p^p+|\mathcal O||c|^p\bigr),
 \quad p'=\frac p{p-1}.
 \end{aligned}
\]
Thus the claimed bound holds after absorbing the fixed factor
$|\mathcal O|$ into the constant.

\smallskip
The mass process is adapted and has moments of every order on finite
intervals; its paths are bounded on each compact interval.
Consequently $c(t)$ is an adapted time-dependent parameter and the
growth bound is integrable. The map $z\mapsto M_tG_0+Cz$ is affine with
bounded linear part on the zero-mean space. These are the variational
hypotheses of \cite{LiuRockner2010}, and so an adapted solution of
\eqref{4.6} is obtained. This argument uses monotonicity only
within the prescribed-mass equation.

\smallskip
Next, we consider \emph{energy, reconstruction, and uniqueness.}
Let
\[
 \Sigma_M'(s,z)=M_sG_0+Cz.
\]
Since the $(i0)$ and $(ij)$ channels are orthogonal in the noise space,
or simply by $(a+b)^2\le2a^2+2b^2$,
\[
 \|\Sigma_M'(s,z)\|_{\mathcal L_2(\mathcal U',H_0)}^2
 \le2|M_s|^2\|G_0\|_{\mathcal L_2}^2
    +2\|C\|_{\mathcal L(H_0,\mathcal L_2)}^2\|z\|_{H_0}^2.
\]
Applying It\^o's formula to $\|z\|_{H_0}^2$ and using
\eqref{4.7} gives
\[
 \begin{aligned}
 \|z(t)\|_{H_0}^2
 +2^{2-m}\int_0^t\|z(s)\|_p^p\dd s
 &\le \|z_0\|_{H_0}^2
 +\int_0^t\|\Sigma_M'(s,z(s))\|_{\mathcal L_2}^2\dd s\\
 &\quad
 +2\int_0^t(z(s),\Sigma_M'(s,z(s))\dd W'_s)_{H_0}.
 \end{aligned}
\]
Define the stopping time $\tau_R=\inf\{s:\|z(s)\|_{H_0}\ge R\}\wedge T$.  The BDG and Young
inequalities imply
\[
 \begin{aligned}
 &\E\sup_{\theta\le t\wedge\tau_R}
 \left|2\int_0^\theta
 (z,\Sigma_M'(s,z)\dd W'_s)_{H_0}\right|\\
 &\le
 \frac12\E\sup_{\theta\le t\wedge\tau_R}\|z(\theta)\|_{H_0}^2
 +C\E\int_0^{t\wedge\tau_R}
     \bigl(|M_s|^2+\|z(s)\|_{H_0}^2\bigr)\dd s .
 \end{aligned}
\]
Therefore Gronwall's inequality gives, with a constant independent of
$R$,
\[
 \E\sup_{s\le T\wedge\tau_R}\|z(s)\|_{H_0}^2
 +\E\int_0^{T\wedge\tau_R}\|z(s)\|_p^p\dd s
 \le C_T\left(\|z_0\|_{H_0}^2+
      \E\int_0^T|M_s|^2\dd s\right).
\]
Letting $R\uparrow\infty$ proves the same estimate without stopping.
Finally,
\[
 \|u(s)\|_p^p=\|z(s)+M_s/|\mathcal O|\|_p^p
 \le2^{p-1}\left(\|z(s)\|_p^p+
      |\mathcal O|^{1-p}|M_s|^p\right),
\]
and the geometric Brownian motion $M$ has moments of every positive
order on finite intervals.  This yields the asserted energy bound for
$u$. Uniqueness follows first for the scalar mass equation and
then for the zero-mean equation with the same prescribed mass, where \eqref{4.7} is used again. The two coordinates reconstruct a solution
of \eqref{4.4} with the boundary condition in
\eqref{4.1}.
If $M_0\ne0$, the explicit exponential in \eqref{4.5} never
vanishes. Since $u(t)=0$ would force $M_t=0$, non-hitting follows.
\end{proof}

\begin{theorem}[Neumann non-hitting from nonzero initial mass]
\label{thm4.4}
If $M(u_0)\ne0$, the solution in Proposition~\ref{prop4.3}
satisfies
$\Pp(u(t)\ne0\text{ for every finite }t\ge0)=1$.
\end{theorem}

\begin{proof}
For Brownian paths in the usual full-probability continuity event,
$\exp\{\sigma\beta_0(t)-\sigma^2t/2\}$ is strictly positive for
every finite $t$. Formula \eqref{4.5} gives $M_t\ne0$ at
all such times. The zero state has mass zero, so it cannot occur.
\end{proof}

Theorem~\ref{thm4.4} covers precisely the class of
nonzero-mass initial states. The mass equation starts at zero and stays
there if $M(u_0)=0$. In that class, the presence of $K_{ij}$ modes and
their nonlocal vacuum ellipticity property do not by themselves provide non-hitting.

\subsection{Partial Malliavin calculus with the mass fixed}
The full Malliavin derivative with respect to $\beta_0$ changes the mass and
creates an indefinite cross term in the zero-mean energy.  The support
argument does not need this direction.  We therefore make the conditional
probability space explicit and differentiate only in the independent
zero-mean noise.

\smallskip
We introduce the Product Wiener space. Because $\beta_0$ and $W'$ are independent, we may define the driving
probability space as
\[
 (\Omega,\mathcal F,\Pp)
 =(\Omega_0\times\Omega',\mathcal F_0\otimes\mathcal F',
   \Pp_0\otimes\Pp').
\]
where $\beta_0=\beta_0(\omega_0)$ and $W'=W'(\omega')$. Since the initial datum and coefficients are deterministic and Proposition~\ref{prop4.3} gives pathwise uniqueness, passing to this canonical product realization does not change the law of the solution, see \cite{BDPRbook}.  Let
$\mathcal U'$ be the noise space generated by the $(i0)$ and $(ij)$ modes,
let $\mathfrak H'=L^2(0,T;\mathcal U')$, and denote by $D'$ the Malliavin
derivative on the Wiener space $(\Omega',\mathcal F',\Pp')$.
For $\Pp_0$-almost every $\omega_0$,
\[
 M^{\omega_0}(s)=M_0\exp\{\sigma\beta_0(s,\omega_0)-\tfrac12\sigma^2s\}
\]
is a bounded continuous deterministic function of $s\in[0,T]$.  With
$c^{\omega_0}=M^{\omega_0}/|\mathcal O|$, equation \eqref{4.6} is then a
variational SPDE on $\Omega'$ with deterministic time-dependent mass
parameter.  We write $\E'$ for expectation on $\Omega'$.  All conditional
Malliavin assertions below mean statements on this second Wiener space for
$\Pp_0$-almost every fixed $\omega_0$.

\smallskip
The reason for using $D'$ can also be seen from the energy.  If a derivative
$Y$ has nonzero mass, the positive-mode $H_0$ pairing of the linearized drift
with its zero-mean part contains
\[
 -\int_{\mathcal O}aY^2 \dd x
 +\overline Y\int_{\mathcal O}aY \dd x,
 \qquad a=m|u|^{m-1},
\]
and the second term has no sign.  In contrast, $D'M=0$ and every $D'z$
has zero mass, so the cross term vanishes.

\begin{proposition}[Conditional terminal estimate]
\label{prop4.5}
There exists a set $\Omega_0^*\in\mathcal F_0$ with
$\Pp_0(\Omega_0^*)=1$ such that the following holds.  Fix
$\omega_0\in\Omega_0^*$, $t>0$, and
$\varphi\in C_c^\infty(\mathcal O)$.  Let
\[
F^{\omega_0}_{t,\varphi}=\langle u^{\omega_0}(t),\varphi\rangle.
\]
Then $F^{\omega_0}_{t,\varphi}\in\mathbb D^{1,2}(\Omega')$ and
\begin{equation}\label{4.8}
 \E'\!\left[
 \left(\frac1\delta\int_{t-\delta}^t
  \|D'_rF^{\omega_0}_{t,\varphi}
        -\Theta_{t,\varphi}^{\prime,\omega_0}\|_{\mathcal U'}^2 \dd r
 \right)^q\right]\longrightarrow0
 \qquad(\delta\downarrow0),
\end{equation}
where $q=(m+1)/(2m)$ as in \eqref{3.4} and
\begin{equation}\label{4.9}
 \Theta_{t,\varphi}^{\prime,\omega_0}
 =\Bigl(
   (M_t^{\omega_0}b_i\langle f_i,\varphi\rangle)_i,
   (a_{ij}\ell_j(z^{\omega_0}(t))
             \langle f_i,\varphi\rangle)_{i,j}
   \Bigr).
\end{equation}
Equivalently, \eqref{4.8} is the conditional terminal estimate given the
mass Brownian path, but its formulation above does not rely on choosing a
version of a regular conditional probability.
\end{proposition}

\begin{proof}
Fix $\omega_0\in\Omega_0^*$ and suppress the superscript $\omega_0$.
Thus $M$ and $c=M/|\mathcal O|$ are deterministic, continuous and bounded
on $[0,T]$. The proof is divided into five steps.

\smallskip
\emph{Step 1: conditional state estimates and Galerkin convergence.}
Let $P_n^N$ be the projection onto the first $n$ positive Neumann modes and
truncate $W'$ to its first $n$ coordinates in a fixed enumeration.  Denote
the resulting zero-mean solution by $z_n$, put $u_n=z_n+c$, and write
\[
 \Sigma_M'(s,z)=M_sG_0+Cz,
\]
where $G_0$ and $C$ are extended by zero to the complementary coordinates
of $\mathcal U'$. For the
fixed path $\omega_0$, set
\[
 K_M(\omega_0)=\sup_{0\le s\le T}|M_s^{\omega_0}|^2<\infty .
\]
The Galerkin equation has the form
\[
 \dd z_n+A_{c(s),n}(z_n)\dd s
 =P_n^N\Sigma_M'(s,z_n)\Pi_n\dd W'_s .
\]
By using \eqref{4.7}, contractivity of $P_n^N,\Pi_n$, and boundedness of $G_0,C$, we get
\[
 \begin{aligned}
 2\langle A_{c(s),n}(z_n),z_n\rangle
 &\ge2^{2-m}\|z_n\|_p^p,\\
 \|P_n^N\Sigma_M'(s,z_n)\Pi_n\|_{\mathcal L_2}^2
 &\le C_0K_M(\omega_0)+C_1\|z_n\|_{H_0}^2,
 \end{aligned}
\]
where $C_0,C_1$ are deterministic and independent of $n$.
Consequently, It\^o's formula gives
\[
 \begin{aligned}
 \|z_n(t)\|_{H_0}^2
 +2^{2-m}\int_0^t\|z_n(s)\|_p^p\dd s
 &\le\|P_n^Nz_0\|_{H_0}^2
 +C_0tK_M(\omega_0)
 +C_1\int_0^t\|z_n(s)\|_{H_0}^2\dd s\\
 &\quad+2\int_0^t
 (z_n,P_n^N\Sigma_M'(s,z_n)\Pi_n\dd W'_s)_{H_0}.
 \end{aligned}
\]
Let
\[
\tau_R^n=\inf\{s:\|z_n(s)\|_{H_0}\ge R\}\wedge T.
\]
BDG and Young's inequality give
\[
 \begin{aligned}
 &\E'\sup_{\theta\le t\wedge\tau_R^n}
 \left|2\int_0^\theta
 (z_n,P_n^N\Sigma_M'(s,z_n)\Pi_n\dd W'_s)_{H_0}\right|\\
 &\le
 \frac12\E'\sup_{\theta\le t\wedge\tau_R^n}
       \|z_n(\theta)\|_{H_0}^2
 +C\E'\int_0^{t\wedge\tau_R^n}
       \bigl(K_M(\omega_0)+\|z_n(s)\|_{H_0}^2\bigr)\dd s .
 \end{aligned}
\]
After absorbing the first term and applying Gronwall inequality, we have
\[
 \E'\sup_{s\le T\wedge\tau_R^n}\|z_n(s)\|_{H_0}^2
 +\E'\int_0^{T\wedge\tau_R^n}\|z_n(s)\|_p^p\dd s
 \le C_T\bigl(\|z_0\|_{H_0}^2+TK_M(\omega_0)\bigr).
\]
The right-hand side is independent of $n$ and $R$.  Fatou's lemma and
monotone convergence as $R\uparrow\infty$ prove
\begin{equation}\label{4.10}
 \sup_n\E'\!\left[
   \sup_{s\le T}\|z_n(s)\|_{H_0}^2
   +\int_0^T\|z_n(s)\|_p^p \dd s\right]
 \le C_T(\omega_0)<\infty.
\end{equation}
The quantity required in the terminal-strip estimate is the $L^p$ norm of
$u_n=z_n+c$. The elementary convexity inequality
$|a+b|^p\le2^{p-1}(|a|^p+|b|^p)$ gives pointwise
\[
 |z_n(s,x)+c(s)|^p
 \le2^{p-1}\bigl(|z_n(s,x)|^p+|c(s)|^p\bigr).
\]
Integrating over $\mathcal O$ yields
\begin{equation}\label{N-u-bound}
 \|z_n(s)+c(s)\|_p^p
 \le 2^{p-1}\bigl(\|z_n(s)\|_p^p
       +|\mathcal O| |c(s)|^p\bigr).
\end{equation}
Since the fixed continuous mass path satisfies
$\sup_{s\le T}|c(s)|<\infty$, integrating this inequality in time and
using \eqref{4.10} gives
\[
 \sup_n\E'\int_0^T\|u_n(s)\|_p^p\dd s
 \le2^{p-1}\left[
  \sup_n\E'\int_0^T\|z_n(s)\|_p^p\dd s
  +T|\mathcal O|\sup_{s\le T}|c(s)|^p\right]<\infty.
\]
Moreover, for every fixed $r\ge2$, applying It\^o's formula to
$(1+\|z_n\|_{H_0}^2)^{r/2}$ yields
\begin{equation}\label{N-high-mom}
 \sup_n\E'\sup_{s\le T}(1+\|z_n(s)\|_{H_0}^2)^{r/2}
 \le C_{r,T}(\omega_0)<\infty,
\end{equation}
Indeed, the first derivative of
$\Psi_r(z)=(1+\|z\|_{H_0}^2)^{r/2}$ is a multiple of
$(1+\|z\|^2)^{r/2-1}z$, and its Hessian is bounded by
$C_r(1+\|z\|^2)^{r/2-1}$ in operator norm.  The drift term is
nonpositive because
$\langle A_{c(s),n}(z_n),z_n\rangle\ge0$.  Furthermore,
\[
 \|\Sigma_M'(s,z_n)\|_{\mathcal L_2}^2
 \le C\bigl(K_M(\omega_0)+\|z_n\|_{H_0}^2\bigr).
\]
Consequently both It\^o correction terms are bounded by
$C_{r,\omega_0}(1+\Psi_r(z_n))$.  After using stopping time trick, the BDG
inequality controls the martingale by one half of $\E'\sup_{q\le s}\Psi_r(z_n(q))$ plus the same integral term.
By using the Gronwall inequality and removal of the stopping times, it is easy to get \eqref{N-high-mom}.

\smallskip
The same monotonicity--Minty argument as in Proposition~\ref{prop2.8}, now
on the fixed triple $L^p_0\subset H_0\subset(L^p_0)^*$, gives after passage
to a subsequence
\begin{equation}\label{N-weak}
 z_n\rightharpoonup z
 \quad\hbox{in }L^p(\Omega'\times(0,T)\times\mathcal O)
 \cap L^2(\Omega'\times(0,T);H_0).
\end{equation}
Uniqueness of the fixed-mass equation identifies the limit with the
conditional solution $z$.  Since the limit is unique, every subsequence has
the same weak limit.

\smallskip
\emph{Step 2: finite-dimensional partial Malliavin derivative.}
For the fixed mass path, the Galerkin drift is
$-A_{c(s),n}$, where
\[
 (DA_{c(s),n}(z)y,y)_{H_0}
 =\int_{\mathcal O}m|z+c(s)|^{m-1}y^2 \dd x\ge0.
\]
To verify this identity, differentiating $\Phi(z+\theta y+c(s))$ at $\theta=0$, we have
\[
 \left.\frac{\dd}{\dd\theta}
 \Phi(z+\theta y+c(s))\right|_{\theta=0}
 =\Phi'(z+c(s))y=m|z+c(s)|^{m-1}y.
\]
Pairing the resulting derivative of $A_{c(s),n}$ with the zero-mean Galerkin vector $y$ gives exactly the displayed integral, and  the positive
Neumann projection disappears in this pairing because $y$ lies in its range. Thus $A_{c(s),n}(z)$ is $C^1$, locally Lipschitz and one-sided Lipschitz with constant zero.
The diffusion $P_n^N\Sigma_M'(s,z)$ is affine with bounded state derivative
$P_n^NC$.  The moment bound \eqref{N-high-mom} verifies the $r>2$ integrability
required in \cite[Corollary~3.5]{ImkellerReisSalkeld2019}.  Hence the
finite-dimensional conditional solution is Malliavin differentiable on
$\Omega'$ without any state cutoff.

\smallskip
For $Y_{r,k}^n=D_r^{\prime k}z_n$ and $s\ge r$,
\begin{equation}\label{N-derivative-eq}
 \begin{aligned}
 \dd Y_{r,k}^n
 &=P_n^N\Delta_N(a_nY_{r,k}^n) \dd s
   +P_n^NCY_{r,k}^n \dd W'_s,\\
 Y_{r,k}^n(r)
 &=P_n^N\Sigma_M'(r,z_n(r))e_k,
 \quad a_n=m|z_n+c|^{m-1}.
 \end{aligned}
\end{equation}
Here the notation $CY \dd W'$ means that $C$ acts only in the $(ij)$
coordinates; it is zero in the $(i0)$ coordinates.  Every $Y_{r,k}^n$ is
zero mean.  Consequently
\[
 (P_n^N\Delta_N(a_nY),Y)_{H_0}
 =-\int_{\mathcal O}a_nY^2 \dd x.
\]
For each fixed $k$, It\^o's formula gives
\[
 \begin{aligned}
 \E'\|Y_{r,k}^n(s)\|_{H_0}^2
 +2\E'\int_r^s\!\!\int_{\mathcal O}a_n|Y_{r,k}^n|^2\dd x\dd q
 &=
 \E'\|Y_{r,k}^n(r)\|_{H_0}^2+
 \E'\int_r^s
 \|P_n^NCY_{r,k}^n(q)\|_{\mathcal L_2}^2\dd q ,
 \end{aligned}
\]
because the stochastic integral has zero expectation after the standard
localization.  Since
\[
\|P_n^NCY\|_{\mathcal L_2}^2\le\|C\|^2\|Y\|_{H_0}^2,
\]
Gronwall inequality implies
\[
 \sup_{r\le s\le T}\E'\|Y_{r,k}^n(s)\|_{H_0}^2
 \le e^{\|C\|^2T}\E'\|Y_{r,k}^n(r)\|_{H_0}^2.
\]
Substitution back into the energy identity controls the
space--time integral by the same initial energy. Summing over $k$ and
using
\[
 \sum_k\|Y_{r,k}^n(r)\|_{H_0}^2
 =\sum_k\|P_n^N\Sigma_M'(r,z_n(r))e_k\|_{H_0}^2
 \le\|\Sigma_M'(r,z_n(r))\|_{\mathcal L_2}^2,
\]
we obtain
\begin{equation}\label{N-derivative-energy}
 \begin{aligned}
 &\sup_{r\le s\le T}\sum_k\E'\|Y_{r,k}^n(s)\|_{H_0}^2
 +\sum_k\E'\int_r^T\!\int_{\mathcal O}
       a_n|Y_{r,k}^n|^2 \dd x\dd s\\
 &\le
 C_T\E'\|\Sigma_M'(r,z_n(r))\|_{\mathcal L_2(\mathcal U',H_0)}^2
 \le C_T(\omega_0).
 \end{aligned}
\end{equation}
The last inequality follows from
\[
\|\Sigma_M'(r,z_n(r))\|_{\mathcal L_2}^2
\le C_0K_M(\omega_0)+C_1\|z_n(r)\|_{H_0}^2
\]
and \eqref{4.10}.  Thus the constant is uniform in $n$ and $r$.
The constants in \eqref{4.10}--\eqref{N-derivative-energy} may depend on the fixed mass
path, but never on $n$ or on the terminal-window length $\delta$.

\smallskip
\emph{Step 3: uniform terminal remainder estimate.}
Set
\[
 \overline\varphi=|\mathcal O|^{-1}\int_{\mathcal O}\varphi \dd x,
 \qquad \psi=\varphi-\overline\varphi,
 \qquad \psi_n=P_n^N\psi.
\]
By the Neumann part of Lemma~\ref{lem2.7},
$\sup_n\|\Delta\psi_n\|_\infty<\infty$ and
$\psi_n\to\psi$ in $H^1$.  For
$F_n(t)=\langle z_n(t),\psi_n\rangle+M_t\overline\varphi$, define
\[
 b_n(r)=\Pi_n\Sigma_M'(r,z_n(r))^*\psi_n.
\]
The derivative equation gives
\[
 D'_rF_n(t)-b_n(r)=R_r^n+S_r^n,
\]
with the drift and stochastic remainders defined exactly as in
\eqref{3.6}, with $a_n$, $\psi_n$ and the zero-mean noise replacing their
Dirichlet counterparts.  Let
\[
 A'_{n,\delta}=\frac1\delta\int_{t-\delta}^t
   \sum_k\int_r^t\!\int_{\mathcal O}a_n|Y_{r,k}^n|^2 \dd x\dd s\dd r,
 \qquad
 Z_n=\int_0^T\|u_n(s)\|_p^p \dd s.
\]
For each $r\in[t-\delta,t]$ and each $k$, weighted
Cauchy--Schwarz gives
\[
 |R_{r,k}^n|^2
 \le
 \left(\int_r^t\!\!\int_{\mathcal O}
       a_n|Y_{r,k}^n|^2\dd x\dd s\right)
 \left(\int_r^t\!\!\int_{\mathcal O}
       a_n|\Delta\psi_n|^2\dd x\dd s\right).
\]
Summing in $k$ and enlarging the second time interval to
$[t-\delta,t]$ yields
\[
 \frac1\delta\int_{t-\delta}^t
 \|R_r^n\|_{\mathcal U'}^2\dd r
 \le A'_{n,\delta}C'_{n,\delta},
\]
where
\[
 C'_{n,\delta}
 :=\int_{t-\delta}^t\!\!\int_{\mathcal O}
      a_n|\Delta\psi_n|^2\dd x\dd s.
\]
By Lemma~\ref{lem2.7},
$\sup_n\|\Delta\psi_n\|_\infty<\infty$, and hence
\[
 C'_{n,\delta}
 \le C_\varphi\int_{t-\delta}^t\!\!\int_{\mathcal O}
      |u_n|^{m-1}\dd x\dd s.
\]
Put $\alpha=(m-1)/(m+1)$.  H\"older's inequality on
$[t-\delta,t]\times\mathcal O$ gives
\[
 \begin{aligned}
 \int_{t-\delta}^t\!\!\int_{\mathcal O}|u_n|^{m-1}\le(\delta|\mathcal O|)^{1-\alpha}
 \left(\int_{t-\delta}^t\!\!\int_{\mathcal O}
       |u_n|^{m+1}\dd x\dd s\right)^\alpha\le C\delta^{1-\alpha}Z_n^\alpha.
 \end{aligned}
\]
Therefore
\[
 C'_{n,\delta}\le
 C_\varphi\delta^{1-\alpha}Z_n^\alpha.
\]
By \eqref{N-derivative-energy},
$\sup_{n,\delta}\E'A'_{n,\delta}<\infty$, while \eqref{N-u-bound} gives
\[
\sup_n\E'Z_n<\infty.
\]
With
$q=(m+1)/(2m)$, the same H\"older calculation as in \eqref{3.7} yields
\[
 \sup_n\E'\bigl(A'_{n,\delta}C'_{n,\delta}\bigr)^q
 \le C(\omega_0)\delta^{1/m}.
\]
It\^o isometry for the stochastic remainder and \eqref{N-derivative-energy} give
\[
 \sup_n\E'\left[
 \left(\frac1\delta\int_{t-\delta}^t
       \|S_r^n\|_{\mathcal U'}^2 \dd r\right)^q\right]
 \le C(\omega_0)\delta^q.
\]
Therefore
\begin{equation}\label{N-remainder}
 \sup_n\E'\left[
 \left(\frac1\delta\int_{t-\delta}^t
       \|D'_rF_n(t)-b_n(r)\|_{\mathcal U'}^2 \dd r\right)^q\right]
 \le C(\omega_0)(\delta^{1/m}+\delta^q).
\end{equation}

\smallskip
\emph{Step 4: fixed-time weak convergence and the closed graph of $D'$.}
The weak convergence \eqref{N-weak} alone is not enough for Malliavin
calculus at the prescribed time $t$, so we record the fixed-time step.
For the scalar Galerkin observations, the equation, the uniform flux bound,
BDG and \eqref{4.10}--\eqref{N-u-bound} imply
\[
 \sup_n\E'|F_n(t)-F_n(s)|
 \le C_{\omega_0,\varphi}
       (|t-s|^{1/p}+|t-s|^{1/2}).
\]
The one-sided time-averaging argument used in Step~3 of
Proposition~\ref{prop2.8} therefore gives
\begin{equation}\label{N-fixed-time}
 F_n(t)\rightharpoonup F(t)
 \qquad\hbox{in }L^2(\Omega').
\end{equation}
The estimate \eqref{N-derivative-energy}, paired with $\psi_n$, gives
\[
 \sup_n\E'\int_0^t\|D'_rF_n(t)\|_{\mathcal U'}^2 \dd r<\infty.
\]
Hence, after extraction,
\[
D'F_n(t)\rightharpoonup V \quad {\rm in} \quad
L^2(\Omega';\mathfrak H').
\]
The Malliavin derivative
\[
D':L^2(\Omega')\supset\mathbb D^{1,2}(\Omega')
\to L^2(\Omega';\mathfrak H')
\]
is a closed linear operator
\cite{Nualart2006}; its graph is therefore a closed linear subspace of a
Hilbert space and is weakly closed.  Together with \eqref{N-fixed-time}, this
implies
\[
 F(t)\in\mathbb D^{1,2}(\Omega'),\qquad V=D'F(t).
\]

For fixed $\delta>0$, affine linearity of $\Sigma_M'$ and the strong
convergence of the spectral/noise projections imply
$b_n\rightharpoonup b$ in
$L^2(\Omega'\times(t-\delta,t);\mathcal U')$, where
\[
 b(r)=\Sigma_M'(r,z(r))^*\psi.
\]
Let $\rho=2q=(m+1)/m>1$.  The continuous embedding
\[
L^2(\Omega';L^2(t-\delta,t;\mathcal U'))
\hookrightarrow L^\rho(\Omega';L^2(t-\delta,t;\mathcal U'))
\]
transports weak convergence to the latter reflexive Banach space.  Its norm
is weakly lower semicontinuous, so \eqref{N-remainder} passes to the limit:
\begin{equation}\label{N-limit-rem}
 \E'\left[
 \left(\frac1\delta\int_{t-\delta}^t
       \|D'_rF(t)-b(r)\|_{\mathcal U'}^2 \dd r\right)^q\right]
 \le C(\omega_0)(\delta^{1/m}+\delta^q).
\end{equation}

\smallskip
\emph{Step 5: replacement of the center by its terminal value.}
The fixed-mass solution belongs to
$L^2(\Omega';C([0,T];H_0))$, $M$ is continuous, and
$\Sigma_M'$ is affine and bounded.  Hence
\[
 \omega_{t,\omega_0}(\delta)
 :=\frac1\delta\int_{t-\delta}^t
   \E'\|b(r)-b(t)\|_{\mathcal U'}^2 \dd r
 \longrightarrow0.
\]
Since $0<q<1$, the squared triangle inequality followed by
$(a+b)^q\le a^q+b^q$ gives from \eqref{N-limit-rem}
\[
 \E'\left[
 \left(\frac1\delta\int_{t-\delta}^t
  \|D'_rF(t)-b(r)\|_{\mathcal U'}^2 \dd r\right)^q\right]
 \le C(\omega_0)
   \bigl(\delta^{1/m}+\delta^q+\omega_{t,\omega_0}(\delta)^q\bigr)
 \longrightarrow0.
\]
Finally,
\[
 F(t)=\langle z(t),\varphi-\overline\varphi\rangle
       +M_t\overline\varphi,
 \qquad D'M_t=0,
\]
and $b(t)$ is precisely the vector in \eqref{4.9}.  This proves
\eqref{4.8}.
\end{proof}

\begin{theorem}[Neumann full-domain support for nonzero initial mass]
\label{thm4.6}
Assume \eqref{4.2}--\eqref{4.3} and let
$u_0\in\mathcal H_N$ be deterministic with $M(u_0)\ne0$.
For every fixed $t>0$ and every nonzero
$\varphi\in C_c^\infty(\mathcal O)$, the scalar observation
$\langle u(t),\varphi\rangle$ has an absolutely continuous law.
Moreover,
\begin{equation}\label{4.11}
 \Pp\bigl(\supp_{\mathcal O}u(t,\cdot)=\overline{\mathcal O}\bigr)=1.
\end{equation}
No separate non-hitting or Malliavin assumption is imposed.
\end{theorem}

\begin{proof}
By Theorem~\ref{thm4.4}, $M_t\ne0$ for every finite $t$ on a
$\Pp_0$-full set of mass paths.  Fix such a path in the full set
$\Omega_0^*$ of Proposition~\ref{prop4.5}.  Since every $b_i$ is positive
and the span of $(f_i)$ is dense in $L^2_0(\mathcal O)$, a nonzero
$\varphi\in C_c^\infty(\mathcal O)$ satisfies
\[
 \|\Theta_{t,\varphi}^{\prime,\omega_0}\|_{\mathcal U'}^2
 \ge (M_t^{\omega_0})^2
      \sum_i b_i^2|\langle f_i,\varphi\rangle|^2>0.
\]
Indeed, if every pairing vanished then the zero-mean part of $\varphi$
would vanish, so $\varphi$ would be constant on the connected domain;
compact support then forces $\varphi=0$.

\smallskip
Now, we apply Lemma~\ref{lem2.13} on the Wiener space
$(\Omega',\mathcal F',\Pp')$, with $D'$ in place of $D$ and with the
terminal estimate of Proposition~\ref{prop4.5}.  The conditional Malliavin
covariance is strictly positive $\Pp'$-almost surely, hence
$F_{t,\varphi}^{\omega_0}$ has a density with respect to Lebesgue measure.
Thus for every Lebesgue-null Borel set $A$,
\[
 \Pp'\bigl(F_{t,\varphi}^{\omega_0}\in A\bigr)=0
 \qquad\hbox{for }\Pp_0\hbox{-a.e. }\omega_0.
\]
Fubini's theorem on the product realization gives
\[
 \Pp(F_{t,\varphi}\in A)
 =\int_{\Omega_0}
   \Pp'\bigl(F_{t,\varphi}^{\omega_0}\in A\bigr)\Pp_0(\dd\omega_0)
 =0,
\]
which proves the unconditional absolute continuity. In
particular $\Pp(F_{t,\varphi}=0)=0$
for every nonzero localized test.
Lemma~\ref{lem2.15} then yields \eqref{4.11}.
\end{proof}

\begin{remark}\label{rem4.7}
(i) The nonzero-mass restriction in Theorem~\ref{thm4.6} is
essential to this proof. The mass formula alone says nothing about
non-hitting for nonzero data of zero mass. In the nonzero-mass class,
the $K_{ij}$ modes are not needed for the conclusion: one may set them
all to zero and retain the same proof. They were included to make the
entire family satisfy nonlocal vacuum ellipticity also at nonzero zero-mass states, as verified
in Proposition~\ref{prop4.2}.

\smallskip
(ii) If also $\sigma=0$, mass is constant and the $K_{i0}$ forcing becomes
the additive coefficient $M_0G_0$ along the solution. For $\sigma\ne0$,
the preceding conditional proof treats its random amplitude. In both
cases the original state-to-noise maps $v\mapsto M(v)b_if_i$ are linear
and spatially nonlocal.
\end{remark}

\begin{remark}\label{rem4.8} (i) The distinction between the Dirichlet and Neumann conclusions lies in
whether the zero state can be excluded. In the Dirichlet variational
space $H_D=(H_0^1)^*$, the constant function $1$ is not an admissible
test function, so the mass argument used for the Neumann problem does
not apply. Under \eqref{3.1} and nonlocal vacuum ellipticity, the Dirichlet argument gives the
terminal estimate but does not establish non-hitting. In the Neumann
product space $\mathcal H_N=\mathbb R\oplus H_0$, the mass satisfies a
closed scalar SDE. Nonzero initial mass therefore rules out the zero
state.

\smallskip
(ii) The Malliavin arguments also use different noise directions. The
Dirichlet estimate involves all driving directions and controls the
terminal quantity $B(u(t))^*\varphi$. In the Neumann case, we condition
on the mass Brownian motion and differentiate only along the
independent zero-mean noise directions. These directions suffice to
obtain strictly positive conditional covariance.
\end{remark}

\section{Numerical Example}
\label{sec:numerical-propagation}
This section gives a one-dimensional numerical illustration of the propagation
mechanisms discussed above.  The experiment is designed to distinguish three
effects: finite propagation for the deterministic porous medium equation,
immediate excitation of an interior vacuum by additive noise, and immediate
excitation of the same vacuum by a genuinely nonlocal multiplicative
coefficient.  The purpose is illustrative rather than a numerical proof of the
Malliavin-support theorems.

\smallskip
The bounded-domain theory is formulated for a bounded connected $C^2$ domain
$\mathcal O=B_1(0)\subset\R^2$.
The purpose is to visualize the same three mechanisms as in the analytical
results: persistence of an interior vacuum for the deterministic porous-medium
flow, immediate vacuum activation by additive noise, and vacuum activation by
a spatially nonlocal multiplicative coefficient.

\smallskip
Let
\[
 m=2,\qquad
 u_0(x)=\left(1-\frac{|x|^2}{0.30^2}\right)_+.
\]
Thus $\supp u_0\subset B_{0.30}(0)$.  We compare
\begin{align*}
 \partial_tu^{\rm D}&=\Delta\Phi(u^{\rm D}),\ \
 \dd u^{\rm A}=\Delta\Phi(u^{\rm A})\dd t
   +\sum_{k=1}^{7}\sigma_k^A\psi_k\dd\beta_k(t),\\
 \dd u^{\rm N}&=\Delta\Phi(u^{\rm N})\dd t
   +\sum_{k=1}^{7}\sigma_k^N\Psi_k(u^{\rm N})\dd\beta_k(t),
\end{align*}
where $\Phi(r)=|r|r$.  For the numerical additive modes we take
\[
 \psi_k(x)=(1-|x|^2)_+
 \frac1{2\pi\ell^2}
 \exp\!\left(-\frac{|x-q_k|^2}{2\ell^2}\right),
 \qquad \ell=0.25,
\]
with shifts
\[
 q_k\in\{0,\ \pm0.4e_1,\ \pm0.4e_2,\ (0.3,0.3),\ (-0.3,-0.3)\},
\]
and weights $\sigma_k^A=0.035k^{-1.1}$.
For the multiplicative model we use the finite-mode nonlocal coefficient
\[
 \Psi_k(v)(x)=(1-|x|^2)_+ (g_{q_k}*v)(x),
 \qquad \sigma_k^N=0.80k^{-1.1}.
\]
The cutoff only enforces compatibility with the boundary in the finite-grid
experiment; the map remains spatially nonlocal in the interior because its
value at $x$ depends on the state over the whole disk.

\begin{table}[htbp]
\centering
\small
\renewcommand{\arraystretch}{1.2}
\begin{tabular}{lccc}
\toprule
Quantity & Deterministic & Additive & Nonlocal multiplicative\\
\midrule
$\mathcal A_2(\Delta t)$
& $0$
& $3.2152\times10^{-6}$
& $6.2839\times10^{-6}$\\
$\mathcal A_2(T)$
& $0$
& $1.4170\times10^{-5}$
& $4.0149\times10^{-5}$\\
active annulus fraction at $\Delta t$
& $0$
& $0.9969$
& $1.0000$\\
active annulus fraction at $T$
& $0$
& $1.0000$
& $1.0000$\\
\bottomrule
\end{tabular}
\caption{Two-dimensional unit-disk computation.  A grid point is called
active when $|u|>10^{-8}$.}
\label{tab:num2d}
\end{table}

\begin{figure}[htbp]
\centering
\includegraphics[width=0.60\textwidth]{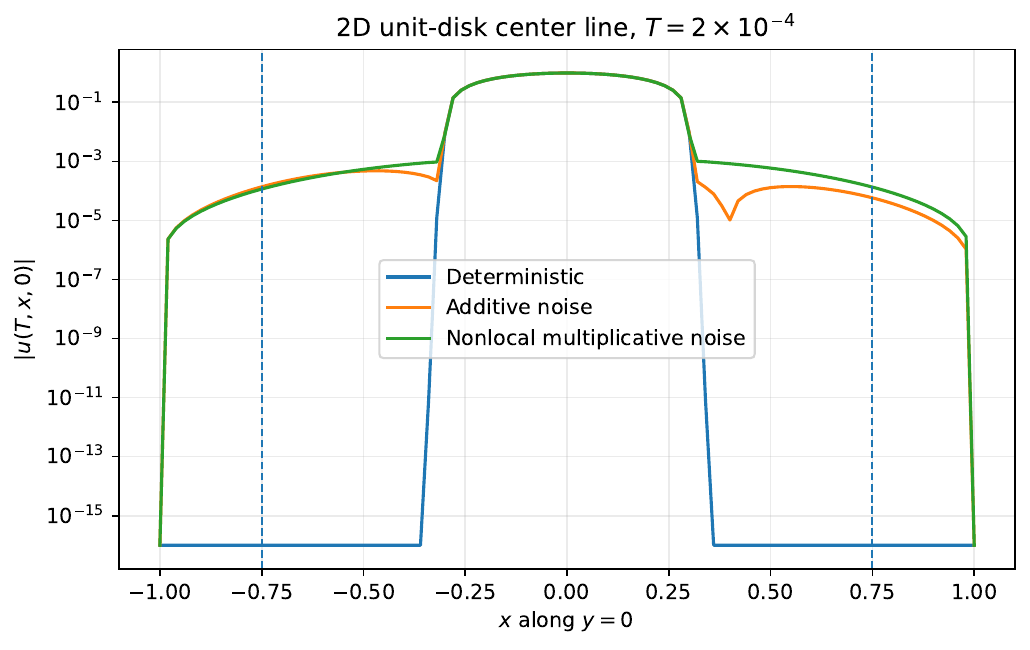}
\caption{Absolute value along the diameter $y=0$ at the final time.}
\label{fig:num2dline}
\end{figure}

\begin{figure}
\centering
		\subfigure[]{
			\includegraphics[width=0.3\linewidth]{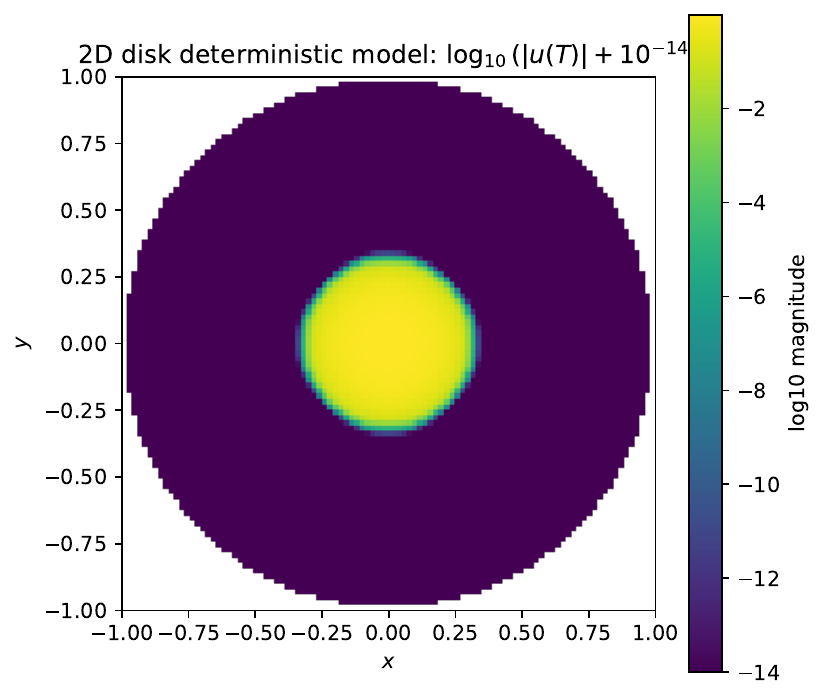}}
\centering
		\subfigure[]{
			\includegraphics[width=0.3\linewidth]{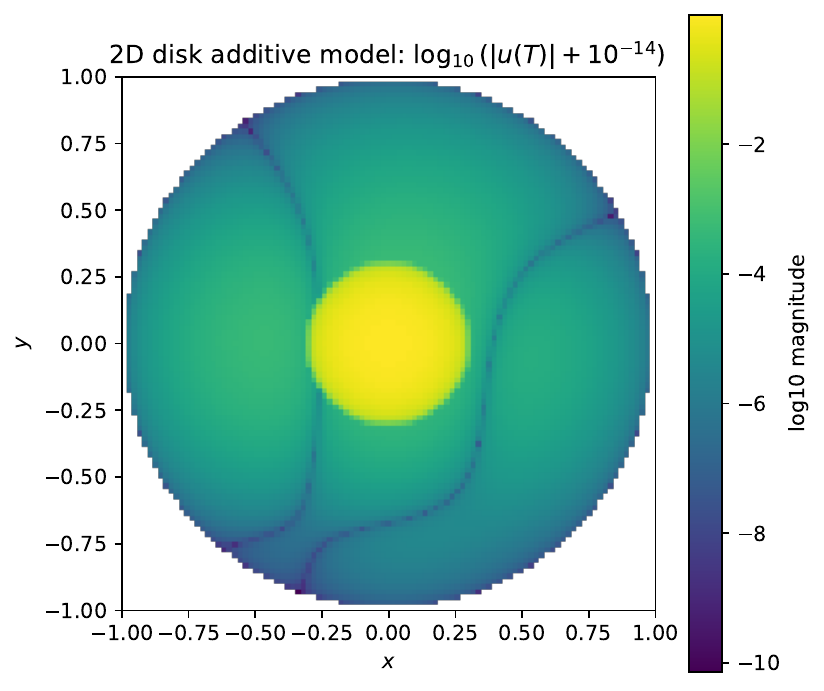}}
\centering
		\subfigure[]{
			\includegraphics[width=0.3\linewidth]{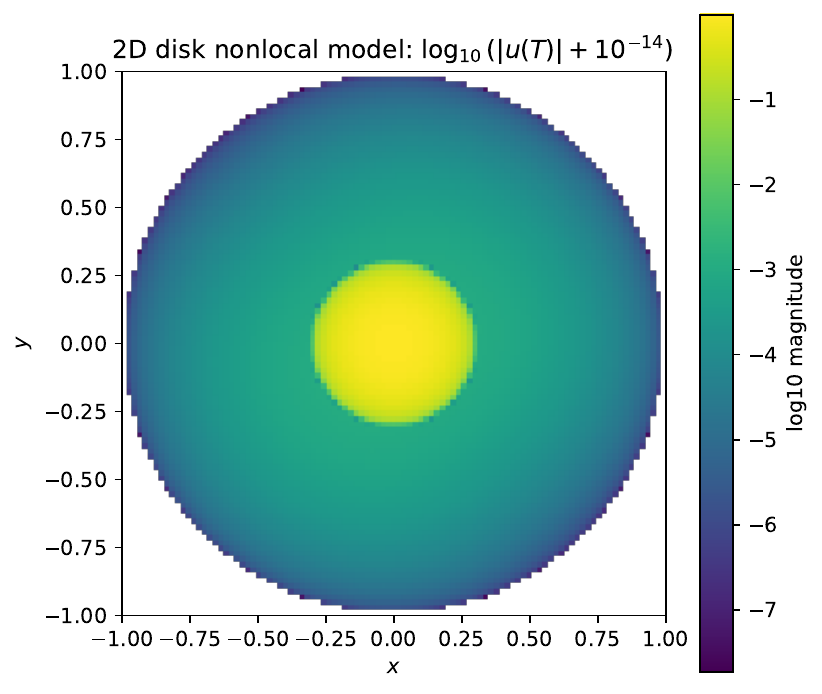}}\vspace{0.5cm}
\caption{(a)--(c) respectively represent the logarithmic magnitude of one realization of deterministic, additive, and
nonlocal multiplicative models.}
\label{2D}
\end{figure}

A Cartesian grid with spacing $\Delta x=0.02$ is masked by the unit disk and
zero values are imposed outside the disk.  We use $\Delta t=10^{-5}$, $T=2\times10^{-4}$,
and Euler--Maruyama for the stochastic increments.  Define the initially empty
remote annulus
\[
 A_{\rm rem}^{(2)}=\{x:0.75\le |x|\le0.90\}
\quad
{\rm and} \quad
 \mathcal A_2(t)=\int_{A_{\rm rem}^{(2)}}|u(t,x)| \dd x.
\]
The numerical values produced by the accompanying Python code are
summarized in Table \ref{tab:num2d}. The deterministic computation retains an exact remote vacuum on the grid over
the displayed time interval: $\mathcal A_2(T)=0$.
In contrast, both stochastic
models create a nonzero far-field signal after the small time.
The deterministic remote annulus remains empty on the displayed
time interval, while both stochastic models activate it immediately, see Figure \ref{tab:num2d}. In order to see it more clearly,  the exact value of $u(t,x,0)$ is given in Figure \ref{fig:num2dline}. Moreover, for
three models, the Logarithmic magnitude of one realization is shown in Figure \ref{2D}.

\begin{figure}[htbp]
\centering
\includegraphics[width=0.60\textwidth]{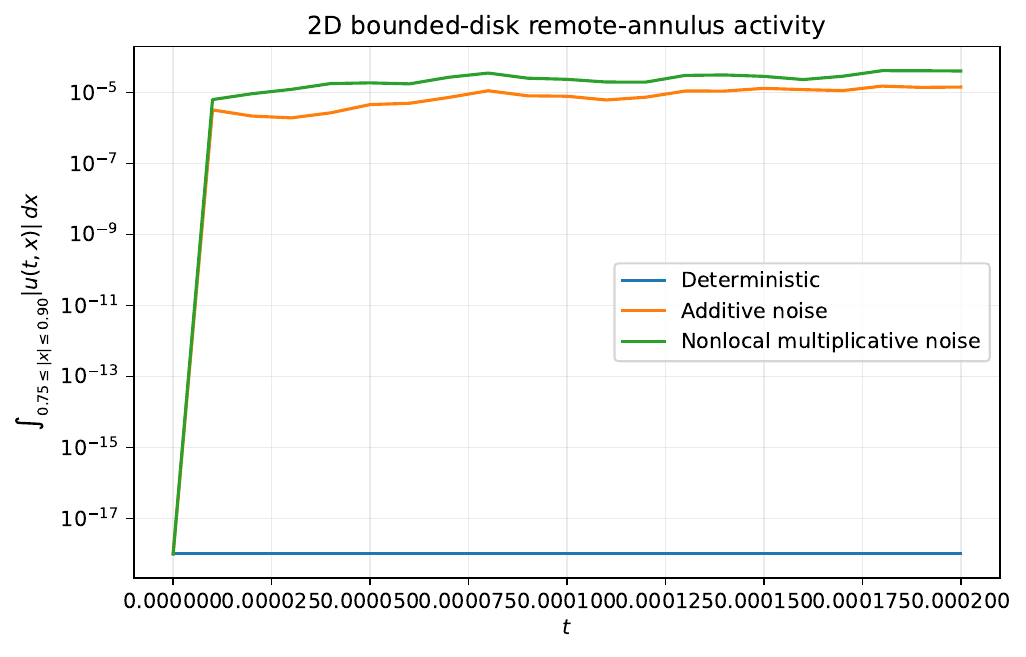}
\caption{Remote-annulus activity in the two-dimensional bounded-domain
experiment.  The deterministic remote annulus remains empty on the displayed
time interval, while both stochastic models activate it immediately.}
\label{fig:num2dactivity}
\end{figure}

\smallskip
Lastly, we give some remarks. The experiment reproduces the qualitative distinction proved analytically.
For the deterministic PME, the diffusion coefficient degenerates at the
vacuum and the interface moves only a finite distance during a short time
interval. Additive noise is independent of the current state, so its
Dirichlet eigenmodes directly inject a signal into a remote vacuum.  The
operator $g_{q_k}*v$ is state dependent but spatially nonlocal.

\smallskip
There are three limitations: First, a spatial grid and a finite noise truncation cannot prove the
continuum statement $\operatorname{supp}_{\mathcal O}u(t)=\overline{\mathcal O}$ almost surely.
They only visualize its mechanism. Second, a support radius in a stochastic
signed solution is threshold dependent.  Third, the
finite nonlocal truncation does not itself satisfy the infinite-dimensional
all-test nonlocal vacuum ellipticity condition.  The analytical Malliavin argument, rather than this
experiment, is what upgrades terminal nondegeneracy to the fixed-time
full-support conclusion.

\section{Conclusion}
We have studied the spatial support of solutions to porous medium
equations on bounded domains under additive and nonlocal multiplicative
noise. For additive noise, we prove that the solution has full
distributional support at each fixed positive time. Since the noise
need not preserve nonnegativity, the argument uses test functions
rather than the positivity set of the solution.

\smallskip
For the Dirichlet problem with nonlocal multiplicative noise, we prove
that the solution at each fixed positive time is either identically
zero or has full distributional support. A non-hitting condition
excludes the zero state and gives full support. In the Neumann model,
the total mass satisfies a closed scalar equation; nonzero initial
mass excludes the zero state without a separate non-hitting assumption. These results show how the support conclusion depends on the noise
and the boundary condition. Whether a similar conclusion holds for
porous medium equations with reflection \cite{RWZ2013} remains an
interesting question.

\section*{Use of AI tools declaration}
During the preparation of this work, the authors used GPT-5.6 Sol to assist with language, the python code and the organization of this paper.

\section*{Acknowledgments}

 This work was supported by the National Natural Science Foundation of China grant 12171247.

\section*{Conflict of interest}

The authors declare there is no conflict of interest.

\section*{Data Availability Statement}

The numerical simulation data are available from the corresponding author upon reasonable request.

\end{document}